\documentclass{article}
\usepackage{todonotes}
\usepackage{geometry}
\usepackage{amsmath}
\usepackage{amssymb}
\usepackage{amsthm}
\usepackage{bm}
\usepackage{graphicx}
\usepackage{color}
\usepackage{tabu}
\usepackage{booktabs}
\usepackage{caption}
\usepackage{authblk}
\usepackage{subfigure}

\usepackage{multirow}

\newtheorem{rem}{Remark}
\newtheorem{defi}{Definition}

\newcommand{\eps}{\varepsilon}
\newcommand{\bx}{\bm{x}}

\newcommand{\bv}{\bm{v}}

\newcommand{\bn}{\bm{n}}

\newcommand*\diff{\mathop{}\!\mathrm{d}}

\begin{document}
\title{Asymptotic-Preserving Neural Networks for Multiscale Time-Dependent Linear Transport Equations}

\author[1, 2]{Shi Jin}
\author[1, 2, 3]{Zheng Ma}
\author[1, \thanks{Corresponding author: wukekever@sjtu.edu.cn}]{Keke Wu}

\affil[1]{School of Mathematical Sciences, Shanghai Jiao Tong University, Shanghai,
    200240, P. R. China.}
\affil[2]{ Institute of Natural Sciences, MOE-LSC,
    Shanghai Jiao Tong University, Shanghai, 200240, P. R. China.}
\affil[3]{Qing Yuan Research Institute,
    Shanghai Jiao Tong University, Shanghai, 200240, China}

\date{\today}

\maketitle
% \tableofcontents

\begin{abstract}
    In this paper we develop a neural network for the numerical {simulation} of time-dependent linear transport equations with diffusive scaling and uncertainties. The goal of the network is to resolve the computational challenges of curse-of-dimensionality and  multiple scales of the problem. We first show that a standard {Physics-Informed Neural Network (PINN)} fails to capture the multiscale nature of the problem, hence justifies the need to use Asymptotic-Preserving  Neural Networks (APNNs). We show that  not all
    classical AP formulations are fit for the neural network approach.
    We construct a micro-macro decomposition based neural network, and also build in a mass conservation mechanism into the loss function, in order to capture the dynamic and multiscale nature of the solutions. Numerical examples are used to demonstrate the effectiveness of this APNNs.
\end{abstract}

\section{Introduction}

% \todo{1. Review of Multiscale kinetic problems and AP schemes.}
In the category of multiscale modeling, kinetic equations  are bridges between continuum and atomistic models \cite{Cerci}. There are two major challenges in kinetic modeling. First is the curse-of-dimensionality, since kinetic equations describe the evolution of probability density function of large number of particles, thus are defined in the phase space, typically a six-dimensional problem plus time.
When uncertainties are considered the dimension can be much higher \cite{Jin-Pareschi-Book, HuJin-UQ, JLM, poette2019gpc, poette2022numerical}. Another difficulty is the multiscale nature. Kinetic equations usually contain small or multiple space and/or temporal scales, characterized by the {Knudsen} number, which is the dimensionless mean free path or time. Often multiscale computations involve  the  coupling of  models at different scales by different numerical schemes \cite{weinan2011principles}. When the Knudsen number is
small,  kinetic  equations can often be approximated by macroscopic hydrodynamic or diffusion equations \cite{BGP}. In this regime,
numerical simulations become prohibitively expensive, since one  needs to  numerically resolve the small physical scales.
Asymptotic-preserving (AP) schemes are those that mimic the asymptotic transitions from the kinetic or hyperbolic equations to their macroscopic (hydrodynamic or diffusive) limits in \emph{discrete setting} ~\cite{jin2010asymptotic}.
Since the mid-1990’s, the development of AP schemes for such problems has generated many interests ~\cite{jin2010asymptotic, Degond-AP, HJL-AP}.
The AP strategy has been proved to be a powerful and robust technique to address multiscale problems in many kinetic and hyperbolic problems.
The main advantage of AP schemes is that they are very efficient in the hydrodynamic or diffusive regime, since they do not need to resolve the small physical parameters numerically and yet can still capture the macroscopic behavior governed by the  hydrodynamic or diffusion equations. In fact, carefully constructed AP schemes avoid the difficulty of coupling a microscopic solver with a macroscopic one, as the micro solver automatically becomes a macro solver in the zero Knudsen number limit.
It was proved, in the case of linear transport with a diffusive scaling, the AP scheme converges uniformly with respect to the scaling parameter~\cite{GJL}.

Due to {aforementioned} computational challenges of kinetic modeling and computations,  efficient computational methods that can efficiently deal with both curse-of-dimensionality and multiple scales are highly desirable. Classical tools to address such issues often use Monte-Carlo methods, which have low-order accuracy, and become
exceedingly expensive when the Knudsen number becomes small \cite{Caf-Par, DP-Acta, PR-MC}. In this paper we seek a machine learning based method to tackle these challenges.

% \todo{2. Review of DNN to solve PDEs, including PINN, DeepRitz, etc.}
The idea of using machine learning method, especially with deep neural networks (DNNs) to solve high-dimensional PDEs has been developed rapidly and achieved some success recently in many problems \cite{E2018, raissi2019physics, lu2021learning, li2020fourier}. In all these approaches, the DNNs came down to minimize the loss function, a high dimensional non-convex optimization problem. The choices of losses usually make a difference for the PDEs; see \cite{E2018, raissi2019physics, li2020fourier, beck2020overview, liao2019deep, deepGalerkin2018,  zang2020weak, cai2021least, lyu2020mim} for reviews and references therein. For applications in kinetic type equations see
\cite{HJJL, CLM, LY, DHY, lou2021physics}.
{The machine learning/DNNs approach has several advantages. First,
it could deal with high-dimensional PDEs due to the strong capacity/representativeness of DNNs.}
Second, it is a mesh-free method so it can deal with problems in complex domain and geometry.
    {Third, it is user friendly to implement for equations by treating the residual error of PDEs as loss function.
        One does not need to construct numerical schemes carefully for approximate derivatives through automatic differentiation technique.}

However, the DNN approach also brings some drawbacks compared with {classical numerical} methods.
The two well-known drawbacks are: lack of high accuracy and {long training time} due to  large number of parameters in DNNs and the use of stochastic gradient method during the training step
for  high-dimensional non-convex functions~\cite{XDE}.
Most importantly, as we find in this article,  that \emph{the standard DeepRitz, PINNs, or NeuralOp approaches have difficulties in dealing  with multiscales, especially for the unsteady  problems.}
This is the main issue we want to address in this paper, and our main strategy to tackle the multiscale challenges is to design an  \emph{Asymptotic-Preserving neural networks} (APNNs).

% \todo[inline]{3. Main idea and the definition of APNNs.}
In this paper, we present a framework of \emph{APNNs} methods for linear transport equations that could exhibit diffusive behavior.
The idea is based on the following observations when using PINNs-like framework to solve multiscale PDEs:
\begin{enumerate}
    \item For the time-dependent/unsteady PDEs, PINNs converts it to a minimization problem of the population/empirical risk/loss in the least square formulation \cite{raissi2019physics}, which \emph{can only capture the leading order/single scale--depending on how to design the loss-- behavior.} This is  not only due to poor selection of the  risk/loss and commonly used first-order gradient based optimization algorithms, but also due to the underlying Frequency-principle \cite{xu2019frequency} of DNN fitting function  that  is difficult to capture the small scale (high frequency) part of the solutions.
    \item Unlike traditional AP schemes, the PINNs framework ignores the asymptotic property of the problem thus cannot capture the macroscopic behavior-which correspond to small physical parameters-- efficiently or effectively. On the contrary, our strategy is to define the  loss based on an AP strategy,  which captures the correct asymptotic behavior when the physical scaling parameters become small, namely the loss has the AP property.
    \item Beyong the AP property, in addition, the conservation or other physical constraints needed to be satisfied by the DNN solutions  to get the correct solutions when dealing with small
          scaling parameters. We propose such a conservation mechanism in our loss.
\end{enumerate}
We first define the  \emph{APNNs}:
\begin{defi}
    Assume the solution is parameterized  using DNN in certain way and trained by using a gradient-based method to minimize a \emph{Loss/Risk}. Then we say it is APNN if, as the physical scaling parameter tends to zero, the loss of the microscopic equation  converges to the loss of the macroscopic equation. That is, the loss when viewed as a numerical approximation of the original equation, has the  AP property.
\end{defi}

% \todo[inline]{4. Main contributions}
Our main contributions lie in the following aspects:
\begin{enumerate}
    \item We propose the concept of \emph{Asymptotic-preserving neural networks} (APNNs) which is a class of machine learning methods that can solve efficiently multiscale equations whose scaling parameter may differ in several orders of magnitude. By requiring loss to satisfy certain conditions (including  AP and conservation), these methods can solve the equations accurately, efficiently and use number of neurons and layers independent of  the scaling parameter. Moreover, it can automatically capture the limiting macroscopic behavior as the scaling parameter tends to zero.
    \item A novel deep neural {network} with structure is constructed for the linear transport equations that exhibit multiscale or diffusive behavior. By using the micro-macro decomposition formulation in the loss and building in the mass-conservation structure, we obtain a machine learning approach to solve the equation which is \emph{APNN}.
\end{enumerate}

% \todo[inline]{5. Structure of the  paper.}
This paper is organized as follows. In Section 2, a brief introduction of the linear transport equation and its diffusion limit is given. The AP micro-macro decomposition is presented in Section 3. Section 4 is our main part where APNNs is introduced with a detailed illustration. Specially, the loss function we choose is based on the micro-macro decomposition and we present our deep neural network that preserves both the asymptotic diffusion limit as well as mass conservation. Numerical problems for both multiscale and high-dimensional uncertainties are given in Section 5, which are solved by APNNs  and compared with PINNs. The paper is concluded in Section 6.

\section{The transport equation and its diffusion limit}
Consider the multidimensional transport equation with the diffusive scaling.
Let $f(t, \bm{x}, \bm{v})$ be the density distribution for particles at space point $\bm{x} \in \mathcal{D} \subset \mathbb{R}^d$, time $t \in [0, T] := \mathcal{T} \subset \mathbb{R}$, and traveling in direction $\bm{v} \in \Omega \subset \mathbb{R}^d$, with $\int_{\Omega} \diff \bm{v} = S$. Here $\Omega$ is symmetric in $\bm{v}$, meaning that $\int_\Omega g(\bm{v}) \, \diff \bm{v} = 0$ for any function $g$ odd in $\bm{v}$.
Then $f$ solves the linear transport equation
\begin{equation}\label{eq:linear-trasport}
    \varepsilon \partial_t f + \bm{v} \cdot \nabla_x f = \frac{1}{\varepsilon} \mathcal{L} f + \varepsilon Q,
\end{equation}
with the collision operator
\begin{equation}\label{eq:linear-trasport operator}
    \mathcal{L} f =  \frac{\sigma_S}{S} \int_\Omega f \, \diff{\bm{v}'}  - \sigma f,
\end{equation}
where $\sigma = \sigma(\bx)$ is the total transport coefficient, $\sigma_S = \sigma_S(\bx)$ is the scattering coefficient, {$Q = Q(\bx)$ is the source term, $\varepsilon > 0$ is the dimensionless mean free path (Knudsen number).} Notice the time variable has been rescaled to the long-time scale $O(1/\varepsilon)$. Here,
\begin{equation}
    \sigma_S = \sigma - \varepsilon^2 \sigma_A,
\end{equation}
where $\sigma_A = \sigma_A(\bx)$ is the absorption coefficient. Such an equation arises in neutron transport, radiative transfer \cite{Chand}, wave propagation in random media \cite{RPK}, etc. In all these applications, the scaling gives rise to a diffusion equation as $\varepsilon \to 0$ \cite{BSS}, which is
\begin{equation}\label{eq:diffusion-limit}
    \partial_t \rho = D\nabla_{\bx} \cdot \left ( \frac{1}{\sigma} \nabla_{\bx} \rho \right ) - \sigma_A \rho + Q\,,
\end{equation}
where $\rho=\int_\Omega f(\bm{v}) \, \diff \bm{v}$. For different collision operators, the diffusion coefficient $D$ may be different. For examples, $D = 1/3$ in one-dimensional slab geometry and $D = 1/2$ when $\Omega$ is a unit sphere in two dimension. This diffusion approximation can be obtained by using the Hilbert expansion with the ansatz
\begin{equation}
    f = f^{(0)} + \varepsilon f^{(1)} + \varepsilon^2 f^{(2)} + \cdots
\end{equation}
and balance the terms by the order in $\varepsilon$. The diffusion approximation can also be derived using the so-called even- and odd-parities \cite{jin2000uniformly}, or micro-macro decomposition technique \cite{lemou2008new} and various Asymptotic-preserving schemes can be constructed based on these methods \cite{jin2010asymptotic}.

For example in the $1$-d case, by splitting equation~\eqref{eq:linear-trasport} and define {even}- and odd-parities as
\begin{equation}\label{even-odd-parities}
    \begin{aligned}
        r(t, x, v) & = \frac{1}{2}[f(t, x, v) + f(t, x, -v)],            \\
        j(t, x, v) & = \frac{1}{2\varepsilon}[f(t, x, v) - f(t, x, -v)],
    \end{aligned}
\end{equation}
one can obtain
\begin{equation}\label{eq:even-odd-parities}
    \begin{aligned}
         & \partial_t r + v\partial_x j                       = \frac{\sigma_S}{\varepsilon^2}(\rho - r) - \sigma_A r + Q, \\
         & \partial_t j + \frac{v}{\varepsilon^2}\partial_x r = -\frac{\sigma_S}{\varepsilon^2} j - \sigma_A j.
    \end{aligned}
\end{equation}
As $\varepsilon \to 0$, the first equation gives $r = \rho$ and the second  gives  $v\partial_x r = -\sigma_S j$. {Substituting} back to the first equation will result
\begin{equation}
    \partial_t \rho = -v^2\partial_x (\frac{1}{\sigma_S}\partial_x \rho) - \sigma_A \rho + Q,
\end{equation}
which after integrating over $v$ yields to the diffusion limit equation~\eqref{eq:diffusion-limit}.

One thing we would like to mention here is that not all classical AP formulations, based on which AP schemes have been constructed,  can be used to  construct the correct APNN method with the desired AP property. See Remark \ref{rmk2}. This is an interesting phenomenon we want to emphasis in this paper. As will be seen in the sequel, the micro-macro decomposition based formulation will help us to construct the desired APNN.

\section{The micro-macro decomposition}

In this section we describe the problem setup and then introduce the micro-macro decomposition. Consider the linear transport equation with initial and boundary conditions over a bounded domain $ \mathcal{T} \times \mathcal{D} \times \Omega$:
\begin{equation}\label{eq:linear-trasport-eqn}
    \left\{
    \begin{array}{ll}
        \varepsilon \partial_t f + \bm{v} \cdot \nabla_x f = \frac{1}{\varepsilon} \mathcal{L} f + \varepsilon Q, & (t, \bm{x}, \bm{v}) \in \mathcal{T} \times \mathcal{D} \times \Omega,\,          \\
        \mathcal{B}f = F_{\text{B}},                                                                              & (t, \bm{x}, \bm{v}) \in \mathcal{T} \times \partial \mathcal{D} \times \Omega,\, \\
        \mathcal{I}f = f_{0},                                                                                     & (t, \bm{x}, \bm{v}) \in \{t = 0\} \times \mathcal{D} \times \Omega,
    \end{array}
    \right.
\end{equation}
where $F_{\text{B}}, f_{0}$ are given functions; $\partial \mathcal{D}$ is the boundary of $\mathcal{D}$, and $\mathcal{B}, \mathcal{I}$ are boundary and initial operators, respectively. More precisely we consider the $1$-d and $2$-d cases in this paper.

\paragraph{One-dimensional case}

Consider the one-dimensional transport equation  in slab geometry
\begin{equation}
    \varepsilon \partial_t f + v\partial_x f = \frac{1}{\varepsilon}\left ( \frac{\sigma_S}{2} \int_{-1}^1 f \, \diff v' - \sigma f \right ) + \varepsilon Q, \quad x_L < x < x_R, \quad -1 \leq v \leq 1,
\end{equation}
with in-flow boundary conditions as,
\begin{equation}
    \begin{aligned}
        f(t, x_L, v) & = F_L(v) \quad \text{for} \quad v > 0, \\
        f(t, x_R, v) & = F_R(v) \quad \text{for} \quad v < 0,
    \end{aligned}
\end{equation}
or  periodic boundary condition
\begin{equation}
    f(t, x_L, v) = f(t, x_R, v).
\end{equation}
The initial function is given as a function of $x$ and $v$
\begin{equation}
    f(0, x, v) = f_0(x, v).
\end{equation}

\paragraph{Two-dimensional case}

The two-dimensional case is very similar except the velocity/angular variables are constrained in the unit circle
\begin{equation}
    \varepsilon \partial_t f + \bv\cdot\partial_{\bx} f = \frac{1}{\varepsilon}\left ( \frac{\sigma_S}{2\pi} \int_{|\bv|=1} f \, \diff \bv' - \sigma f \right ) + \varepsilon Q, \quad \bx \in \Gamma \subset \mathbb{R}^2, \quad |\bv| = 1,
\end{equation}
with $\bv = (\xi, \eta), -1\ll\xi,\eta \ll 1, \xi^2 + \eta^2 = 1$. The in-flow boundary condition is,
\begin{equation}
    f(t, \bx, \bv) = F_B(\bx, \bv) \quad \text{for} \quad \bn\cdot\bv < 0, \quad \bx\in\partial \Gamma,
\end{equation}
where $\bn$ is the outer normal of the boundary. The initial condition is
\begin{equation}
    f(0, \bx, \bv) = f_0(\bx, \bv).
\end{equation}

%\subsection{The micro-macro decomposition}
In this section, we describe the micro-macro decomposition formulation for the linear transport equation. It is a useful mechanism to build AP schemes \cite{lemou2008new}. The idea begins with the decomposition of $f$ into the equilibrium part $\rho$ and the non-equilibrium part $g$:
\begin{equation}\label{eq:decomposition}
    f = \rho + \varepsilon g,
\end{equation}
where
\begin{equation}\label{constrain-rho}
    \rho = \frac{1}{S}\int_\Omega f \, \diff{\bm{v}'}.
\end{equation}
The non-equilibrium part $g$ clearly satisfies $\left \langle g \right \rangle = 0$, where
\begin{equation}\label{constrain-g}
    \left \langle g \right \rangle := \frac{1}{S} \int_\Omega g \, \diff{\bm{v}'} = 0.
\end{equation}
Applying equation (\ref{eq:decomposition}) in (\ref{eq:linear-trasport}) gives
\begin{equation}\label{eq:linear-trasport-2}
    \varepsilon \partial_t \rho + \varepsilon^2  \partial_t g + \bm{v} \cdot  \nabla_{\bm{x}} \rho + \varepsilon \bm{v} \cdot \nabla_{\bm{x}} g = \mathcal{L} g  + \varepsilon Q.
\end{equation}
Integrating this equation with respect to $v$ one obtains the following continuity equation:
\begin{equation}\label{eq:macro}
    \partial_t \rho + \nabla_{\bm{x}} \cdot \left \langle  \bm{v}g \right \rangle = Q.
\end{equation}
Define operator $\Pi: \Pi(\cdot)(\bm{v}) =  \left \langle \cdot \right \rangle$ and $I$ is the identity operator. Then an evolution equation on $g$ is found by applying the orthogonal projection $I - \Pi$ to equation (\ref{eq:linear-trasport-2}):
\begin{equation}\label{eq:micro}
    \varepsilon^2 \partial_t g + \varepsilon (I - \Pi)(\bm{v} \cdot \nabla_{\bm{x}} g) + \bm{v} \cdot  \nabla_{\bm{x}} \rho = \mathcal{L} g + (I - \Pi) \varepsilon Q.
\end{equation}
Finally, \eqref{eq:macro} and \eqref{eq:micro} together with the constraint~\eqref{constrain-rho} and \eqref{constrain-g} constitute the micro-macro formulation of (\ref{eq:linear-trasport}):
\begin{equation}\label{eq:micro-macro}
    \left\{
    \begin{aligned}
         & \partial_t \rho + \nabla_{\bm{x}} \cdot \left \langle  {\bm{v}}g \right \rangle = Q,                                                                               \\
         & \varepsilon^2 \partial_t g + \varepsilon (I - \Pi)(\bm{v} \cdot \nabla_{\bm{x}} g) + \bm{v} \cdot  \nabla_{\bm{x}} \rho = \mathcal{L} g + (I - \Pi) \varepsilon Q, \\
         & \langle g \rangle = 0.
    \end{aligned}
    \right.
\end{equation}
or the last equation can be replaced by enforcing it only at initial time
\begin{equation}\label{eq:micro-macro-1}
    \left\{
    \begin{aligned}
         & \partial_t \rho + \nabla_{\bm{x}} \cdot \left \langle  {\bm{v}}g \right \rangle = Q,                                                                               \\
         & \varepsilon^2 \partial_t g + \varepsilon (I - \Pi)(\bm{v} \cdot \nabla_{\bm{x}} g) + \bm{v} \cdot  \nabla_{\bm{x}} \rho = \mathcal{L} g + (I - \Pi) \varepsilon Q, \\
         & \langle g \rangle (0, x) = 0.
    \end{aligned}
    \right.
\end{equation}
This is due to the conservation of $\langle g \rangle$ with respect to $t$ and can be seen easily by integrating the second equation with respect to $v$,
\begin{equation}
    \partial_t \langle g \rangle = 0.
\end{equation}

When $\varepsilon \to 0$, the above system formally approaches
\begin{equation}\label{eq:AP-limit}
    \left\{
    \begin{array}{ll}
        \partial_t \rho + \nabla_{\bm{x}} \cdot \left \langle  {\bm{v}}g \right \rangle = Q, \\
        \bm{v} \cdot  \nabla_{\bm{x}} \rho = \mathcal{L} g .
    \end{array}
    \right.
\end{equation}
The second equation yields
    {
        \begin{equation}\label{q-eqn}
            g=\mathcal{L}^{-1}(\bm{v} \cdot  \nabla_{\bm{x}} \rho),
        \end{equation}}
which, when plugging into the first equation and integrating over $v$, gives the diffusion equation~\eqref{eq:diffusion-limit}.

\begin{rem} In (\ref{eq:micro-macro}) or {(\ref{eq:micro-macro-1})}
    we singled out the condition $\langle g \rangle = 0$. This is not
    necessary in constructing classical AP schemes \cite{lemou2008new}.
    However, for DNN this condition is usually not automatically satisfied and one needs to impose this condition in the loss,
    as shown in the next section.
\end{rem}
\section{Asymptotic-Preserving Neural Networks}

%For the classical Asymptotic-preserving schemes, we focus on the %properties of numerical discretizations as the scaling parameter %changes.
%AP schemes mean to achieve the goal of perserving the numerical %discretizations of the limiting equations as well as uniform %stability.
%Though the design of AP schemes requires understandings of the %specific equation itself and is often case by case, there are some %general ideas that one can follow to construct the schemes %especially kinetic equations.

%However, when we design an APNN method, things change dramatically.
%Using DNNs to solve PDEs follows totally different strategies.
%In this paper, we use the PINNs framework to present the idea of %building an APNN method, however, it is not restricted to only %PINNs and can be generalized to DeepRitz, DeepONet, etc.

\begin{figure}[ht]
    \centering
    \includegraphics[width=0.45\textwidth]{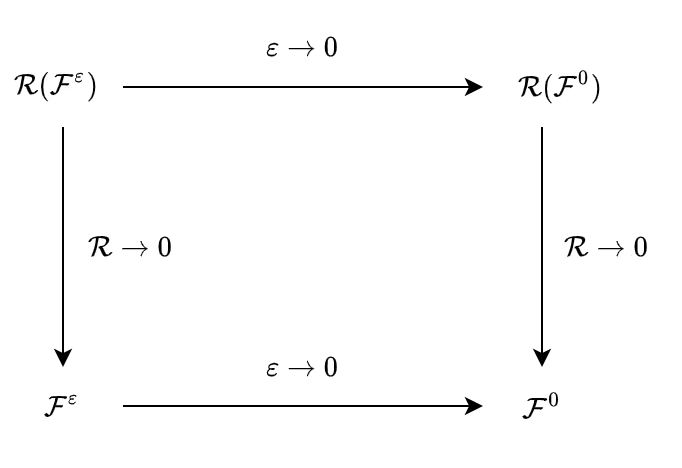}
    \caption{Illustration of APNNs. $\mathcal{F^{\varepsilon}}$ is the microscopic equation that depends on the small scale parameter $\varepsilon$ and $\mathcal{F}^{0}$ is its macroscopic limit as $\varepsilon \to 0$, which is independent of $\varepsilon$. The latent solution of $\mathcal{F^{\varepsilon}}$ is approximated by neural networks with its measure denoted by $\mathcal{R}(\mathcal{F^{\varepsilon}})$. The asymptotic limit of $\mathcal{R}(\mathcal{F^{\varepsilon}})$ as $\varepsilon \to 0$, if exists, is denoted by $\mathcal{R}(\mathcal{F}^{0})$. If $loss(\mathcal{F}^{0})$ is a good measure of $\mathcal{F}^{0}$, then it is called asymptotic-preserving (AP).}
    \label{fig:apnns}
\end{figure}

Unlike classical numerical schemes, DNN framework consists of three ingredients: a neural network parametrization of the solution, a population and empirical loss/risk and an optimization algorithm. In terms of proposed APNNs, the key component  is  to design a loss that has the AP property.  The diagram in Fig.~\ref{fig:apnns}  illustrates the idea of APNNs.
%By comparing this with the traditional AP diagram \cite{GJL, %jin2010asymptotic}, we can think the loss/risk $\mathcal{R}$ as the %discretization or schemes used in the traditional numerical schemes %and think $\mathcal{R}\to 0$ line, that is, the training loss tends %to $0$, as the numerical consistency.

%\subsection{Solving the linear transport equation by APNNs}
The procedure of a deep neural network for solving a PDE problem consists of three parts: a neural network structure, a loss, and a method to minimize loss over the parameter space.

In what follows, conventional notations for deep neural networks (DNNs)\footnote{BAAI.2020.\ Suggested Notation for Machine Learning.\ https://github.com/mazhengcn/suggested-notation-for-machine-learning.} are introduced. An $L$-layer feed forward neural network is defined recursively as,
\begin{equation}
    \begin{aligned}
        f_{\theta}^{[0]}(x) & = x,                                                                              \\
        f_{\theta}^{[l]}(x) & = \sigma \circ (W^{[l-1]} f_{\theta}^{[l-1]}(x) + b^{[l-1]}), \, 1 \le l \le L-1, \\
        f_{\theta}(x)       & = f_{\theta}^{[L]}(x) = W^{[L-1]} f_{\theta}^{[L-1]}(x) + b^{[L-1]},
    \end{aligned}
\end{equation}
where $W^{[l]} \in  \mathbb{R}^{m_{l+1}\times m_l}, b^{l}\in  \mathbb{R}^{m_{l+1}}, m_0 = d_{in} = d$ is the input dimension, $m_{L} = d_0$ is the output dimension, $\sigma$ is a scalar function and $"\circ"$ means entry-wise operation. We denote the set of parameters by $\theta$. For simplicity of neural network presentation, we denote the layers by a list, i.e., $[m_0, \cdots, m_L]$.

There are various choices to build the loss, when applying a DNN to solve a given PDE. In this paper we compare several results about different scale of $\eps$ in  equation~(\ref{eq:linear-trasport-eqn})
with Physics-informed Neural Networks (PINNs) and our proposed Asymptotic Preserving Neural Networks (APNNs). PINNs has the loss as the mean-square error or the residual associated to the given PDE. In contrast, APNNs rewrites the original PDE into a system in asymptotic-preserving form and takes their mean-square residual error as the loss.
Boundary and initial conditions are treated as a regularization or penalty term with penalty parameters $\lambda_1, \lambda_2$ into the loss, which are selected for better performance.

For simplicity in the following discussion we set $\sigma = \sigma_S = 1$.

\subsection{The failure of PINNs  to resolve small scales}
For PINNs, one neural network is applied to directly approximate the density function $f(t, x, v)$,
\begin{equation}
    \text{NN}_{\theta}(t, \bx, \bv) \approx f(t, \bx, \bv).
\end{equation}
The inputs of DNN are $(t, \bm{x}, \bm{v})$, i.e., $m_0 = 3, 5$ for $1$-d and $2$-d respectively. The output is a scalar which represents the value of $f$ at $(t, \bx, \bv)$. Since $f$ is always non-negative, we put an exponential function at the last output layer of the DNN. To be precise, we use
\begin{equation}
    f^{\text{NN}}_{\theta}(t, \bx, \bv) := \exp \left(-\tilde{f}^{\text{NN}}_{\theta}(t, \bx, \bv)\right) \approx f(t, \bx, \bv)
\end{equation}
to represent the numerical solution of $f$. Then the least square of the residual of the original transport equation~\eqref{eq:linear-trasport} is used as the target loss function, together with boundary and initial conditions as penalty terms, that is,
\begin{equation}\label{eq:loss-pinn}
    \begin{aligned}
        \mathcal{R}^{\varepsilon}_{\text{PINN}} = & \frac{1}{|\mathcal{T} \times \mathcal{D} \times \Omega|} \int_{\mathcal{T}} \int_{\mathcal{D}} \int_\Omega \left| \varepsilon^2 \partial_t f^{\text{NN}}_{\theta} + \varepsilon \bm{v} \cdot \nabla_x f^{\text{NN}}_{\theta} - \mathcal{L} f^{\text{NN}}_{\theta} - \varepsilon^2 Q \right|^2 \diff{\bm{v}} \diff{\bm{x}} \diff{t} \\
                                                  & +  \frac{\lambda_1}{|\mathcal{T} \times \partial \mathcal{D} \times \Omega|}  \int_{\mathcal{T}} \int_{\partial \mathcal{D}} \int_\Omega |\mathcal{B}f^{\text{NN}}_{\theta} - F_{\text{B}}|^2 \diff{\bm{v}} \diff{\bm{x}} \diff{t}                                                                                                 \\
                                                  & +  \frac{\lambda_2}{|\mathcal{D} \times \Omega|} \int_{\mathcal{D}} \int_\Omega |\mathcal{I}f^{\text{NN}}_{\theta} - f_{0}|^2 \diff{\bm{v}} \diff{\bm{x}},
    \end{aligned}
\end{equation}
where $\lambda_1$ and $\lambda_2$ are the penalty weights  to be tuned. Then a standard stochastic gradient method (SGD) or Adam optimizer is used to find the global minimum of this loss. Notice that in order to approximate each integral in the loss we use randomly sampled batch at every step, to be detailed in section 5.

Now let us check whether this PINN method is AP.  We only need to focus on the first term of~\eqref{eq:loss-pinn}
\begin{equation}
    \mathcal{R}^{\varepsilon}_{\text{PINN, residual}} := \frac{1}{|\mathcal{T} \times \mathcal{D} \times \Omega|} \int_{\mathcal{T}} \int_{\mathcal{D}} \int_\Omega \left| \varepsilon^2 \partial_t f^{\text{NN}}_{\theta} + \varepsilon \bm{v} \cdot \nabla_x f^{\text{NN}}_{\theta} - \mathcal{L} f^{\text{NN}}_{\theta} - \varepsilon^2 Q \right|^2 \diff{\bm{v}} \diff{\bm{x}} \diff{t}.
\end{equation}
Taking $\varepsilon \to 0$, formally this will lead to
\begin{equation}
    \mathcal{R}_{\text{PINN, residual}} := \frac{1}{|\mathcal{T} \times \mathcal{D} \times \Omega|} \int_{\mathcal{T}} \int_{\mathcal{D}} \int_\Omega \left| - \mathcal{L} f^{\text{NN}}_{\theta} \right|^2 \diff{\bm{v}} \diff{\bm{x}} \diff{t},
\end{equation}
which can be viewed as the PINN loss of the equilibrium equation
\begin{equation}
    \mathcal{L} f = 0.
\end{equation}
This shows that when $\varepsilon$ is very small, to the leading order we are solving equation $\mathcal{L} f = 0$ which gives
$f=\rho$. This does not give the desired   diffusion equation~\eqref{eq:diffusion-limit}. This explains why PINN will fail when $\varepsilon$ is small. We also conduct  numerical experiments to verify this claim in Section 5.

\subsection{The micro-macro decomposition based APNN method}
Now we present an APNN method based on the  micro-macro decomposition method. The main idea is to use PINN to solve the micro-macro system~\eqref{eq:micro-macro} or \eqref{eq:micro-macro-1} {instead} of the original equation~\eqref{eq:linear-trasport}.

First we need to use DNN to parametrize two functions $\rho(x,v)$ and $g(t, x, v)$. So here two networks are used. First
\begin{equation}
    \rho^{\text{NN}}_{\theta}(t, \bx) := \exp \left( -\tilde{\rho}^{\text{NN}}_{\theta}(t, \bx)\right) \approx \rho(t, \bx),
\end{equation}
Notice here $\rho$ is non-negative. Second,
\begin{equation}\label{g-net}
    g^{\text{NN}}_{\theta}(t, \bx, \bv) := \tilde{g}^{\text{NN}}_{\theta}(t, \bx, \bv) - \langle \tilde{g}^{\text{NN}}_{\theta} \rangle (t, \bx) \approx g(t, \bx, \bv).
\end{equation}
Here $\tilde{\rho}$ and $\tilde{g}$ are both fully-connected neural {networks}. Notice that by choosing $g^{\text{NN}}_{\theta}(t, \bx, \bv)$ as in ~\eqref{g-net} it will automatically satisfy the constraint~\eqref{constrain-g} as,
\begin{equation}
    \langle g^{\text{NN}}_{\theta} \rangle = \langle \tilde{g}^{\text{NN}}_{\theta} \rangle - \langle \tilde{g}^{\text{NN}}_{\theta} \rangle = 0, \quad \forall \; t, \bx,
\end{equation}
which automatically satisfies the third equation in the micro-macro system~\eqref{eq:micro-macro}.

Then we propose the least square of the residual of the micro-macro system~\eqref{eq:micro-macro} as the APNN loss,
\begin{equation}\label{eq:loss-ap}
    \begin{aligned}
        \mathcal{R}^{\varepsilon}_{\text{APNN}} = & \frac{1}{|\mathcal{T} \times \mathcal{D}|} \int_{\mathcal{T}} \int_{\mathcal{D}} | \partial_t \rho^{\text{NN}}_{\theta} + \nabla_x \cdot \left \langle   \bm{v} g^{\text{NN}}_{\theta} \right \rangle - Q |^2 \diff{\bm{x}}  \diff{t}                                      \\
                                                  & + \frac{1}{|\mathcal{T} \times \mathcal{D} \times \Omega|} \int_{\mathcal{T}} \int_{\mathcal{D}} \int_\Omega | \varepsilon^2 \partial_t g^{\text{NN}}_{\theta}  + \varepsilon (I - \Pi)(\bm{v} \cdot \nabla_x g^{\text{NN}}_{\theta})                                      \\
                                                  & \quad  + \bm{v} \cdot  \nabla_{\bm{x}} \rho^{\text{NN}}_{\theta} - {\mathcal L} g^{\text{NN}}_{\theta}-(I -\Pi) \varepsilon Q|^2 \diff{\bm{v}} \diff{\bm{x}} \diff{t}                                                                                                      \\
                                                  & +  \frac{\lambda_1}{\mathcal{T} \times\partial \mathcal{D} \times \Omega|}  \int_{\mathcal{T}} \int_{\partial \mathcal{D}} \int_\Omega |\mathcal{B}(\rho^{\text{NN}}_{\theta} + \varepsilon g^{\text{NN}}_{\theta}) - F_{\text{B}}|^2 \diff{\bm{v}} \diff{\bm{x}} \diff{t} \\
                                                  & +  \frac{\lambda_2}{|\mathcal{D} \times \Omega|} \int_{\mathcal{D}} \int_\Omega |\mathcal{I}(\rho^{\text{NN}}_{\theta} + \varepsilon g^{\text{NN}}_{\theta}) - f_{0}|^2 \diff{\bm{v}} \diff{\bm{x}}.
    \end{aligned}
\end{equation}

\begin{figure}[ht]
    \centering
    \includegraphics[width=0.8\textwidth]{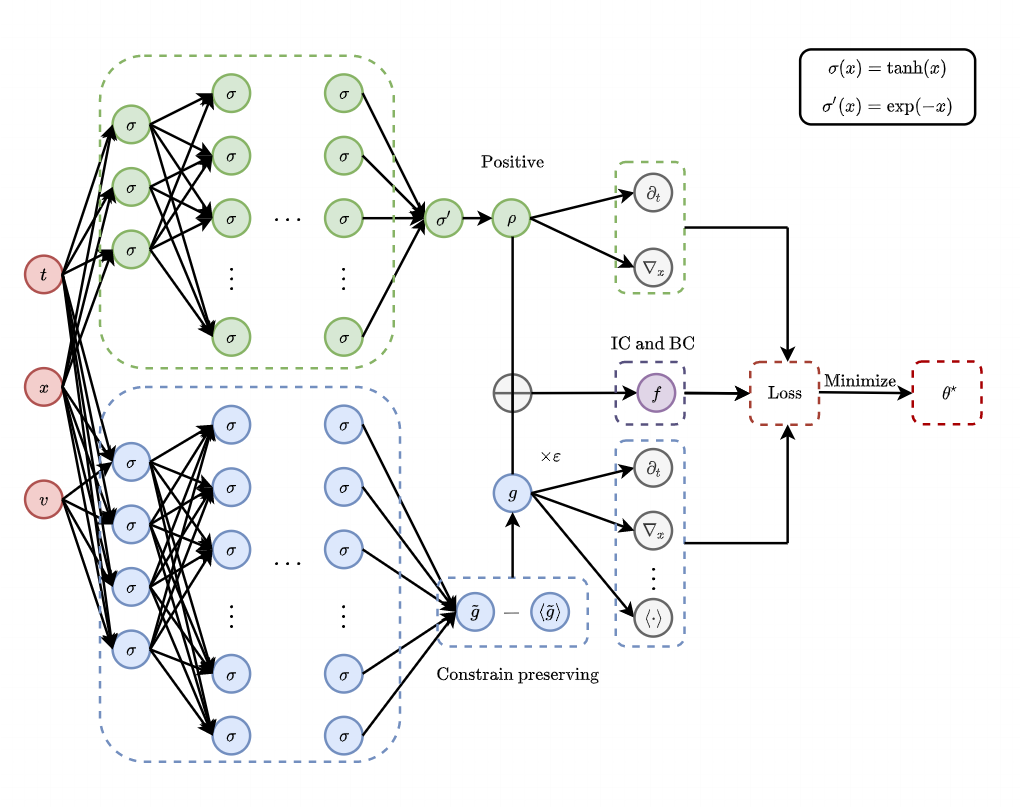}
    \caption{Schematic of APNNs for solving the linear transport equation with initial and boundary conditons.}
    \label{fig:APNNs}
\end{figure}

Now we show formally the AP property of  this loss by considering its behavior for  $\varepsilon$ small. We only need to focus on the first two terms of~\eqref{eq:loss-ap}
{
    \begin{equation}
        \begin{aligned}
            \mathcal{R}^{\varepsilon}_{\text{APNN, residual}} = & \frac{1}{|\mathcal{T} \times \mathcal{D}|} \int_{\mathcal{T}} \int_{\mathcal{D}} | \partial_t \rho^{\text{NN}}_{\theta} + \nabla_x \cdot \left \langle   \bm{v} g^{\text{NN}}_{\theta} \right \rangle - Q |^2 \diff{\bm{x}}  \diff{t}                                                               \\
                                                                & + \frac{1}{|\mathcal{T} \times \mathcal{D} \times \Omega|} \int_{\mathcal{T}} \int_{\mathcal{D}} \int_\Omega \Big | \varepsilon^2 \partial_t g^{\text{NN}}_{\theta} + \varepsilon (I - \Pi)(\bm{v} \cdot \nabla_x g^{\text{NN}}_{\theta}) + \bm{v} \cdot  \nabla_{\bm{x}} \rho^{\text{NN}}_{\theta} \\
                                                                & \qquad \qquad \qquad \qquad \qquad \qquad \qquad - {\mathcal L} g^{\text{NN}}_{\theta}-(I -\Pi) \varepsilon Q \Big|^2 \diff{\bm{v}} \diff{\bm{x}} \diff{t}
        \end{aligned}.
    \end{equation}}
Taking $\varepsilon \to 0$, formally this will lead to
\begin{equation}
    \begin{aligned}
        \mathcal{R}_{\text{APNN, residual}} = & \frac{1}{|\mathcal{T} \times \mathcal{D}|} \int_{\mathcal{T}} \int_{\mathcal{D}} | \partial_t \rho^{\text{NN}}_{\theta} + \nabla_x \cdot \left \langle   \bm{v} g^{\text{NN}}_{\theta} \right \rangle - Q |^2 \diff{\bm{x}}  \diff{t}                           \\
                                              & + \frac{1}{|\mathcal{T} \times \mathcal{D} \times \Omega|} \int_{\mathcal{T}} \int_{\mathcal{D}} \int_\Omega \Big | \bm{v} \cdot  \nabla_{\bm{x}} \rho^{\text{NN}}_{\theta} - {\mathcal L} g^{\text{NN}}_{\theta} \Big|^2 \diff{\bm{v}} \diff{\bm{x}} \diff{t},
    \end{aligned}
\end{equation}
which is the least square loss of equations~\eqref{eq:AP-limit}
\begin{equation}
    \left\{
    \begin{aligned}
         & \partial_t \rho + \nabla_{\bm{x}} \cdot \left \langle  {\bm{v}}g \right \rangle = Q, \\
         & \bm{v} \cdot  \nabla_{\bm{x}} \rho = \mathcal{L} g .
    \end{aligned}
    \right.
\end{equation}
Same as in Section~3.1 the second equation yields $g=\mathcal{L}^{-1}(\bm{v} \cdot  \nabla_{\bm{x}} \rho)$, which, when plugging into the first equation and integrating over $v$, gived the diffusion equation~\eqref{eq:diffusion-limit}.
Hence this proposed method is an APNN method. Finally we put a schematic plot of our method in Figure~\ref{fig:APNNs}.
\begin{rem}
    For the constraint $\left \langle g \right \rangle = 0$, one way is to construct a novel neural network for $g$ such that it  exactly satisfies $\left \langle g \right \rangle = 0$. The other way is to  treat it as a soft constraint with parameter $\lambda_3$, i.e., without using~\eqref{constrain-g}, we use $\hat{g}_{\theta}^{\text{NN}}$ and modifies the loss as
    \begin{equation}\label{eq:soft-constraint-g}
        \mathcal{R}_{\text{APNN, residual}} +  \frac{\lambda_3}{|\mathcal{T} \times \mathcal{D}|} \int_{\mathcal{T}} \int_{\mathcal{D}} | \left \langle  \tilde{g}^{\text{NN}}_{\theta} \right \rangle - 0|^2 \diff{\bm{x}}  \diff{t}.
    \end{equation}
    Our APNNs belongs to the first scenario. Specifically,  in $\mathcal{R}_{\text{apnn}}$, the integrand is not exactly the micro-macro decomposition in (\ref{eq:micro-macro}). Instead we replaced $g$ by $g-\langle{g}\rangle$ since the later always satisfies $\langle{\cdot}\rangle=0$, the desired conservation property!  There are other terms in $\mathcal{R}^{\varepsilon}_{\text{APNN}}$ where such a change will not have any impact since the corresponding operators are invariant under such a transformation.
\end{rem}

\begin{rem}\label{rmk2} The parity formulation (\ref{eq:even-odd-parities}) also allows one to construct AP schemes \cite{jin2000uniformly}. However if one uses it in the loss, it will not be an APNN. To see this, consider the case when $\varepsilon \to 0$. Since the DNN will only pick up the leading term, thus one has, in the $L^2$ sense,
    {
            \begin{equation}\label{parity-limit}
                \begin{aligned}
                     & 0 =\rho - r,                           \\
                     & j = -\frac{v}{\sigma_S}\partial_x r\,.
                \end{aligned}
            \end{equation}}
    These two equations will not lead to the diffusion equation~\eqref{eq:diffusion-limit}.
\end{rem}

\section{Numerical examples}

In this section, in order to verify and compare the performance of PINNs and APNNs, we present both $1$D and $2$D numerical results for several problems chosen from rarefied regimes ($\varepsilon \approx O(1)$) to hydrodynamic (diffusive) regimes ($\varepsilon \to 0$), including examples with high-dimensional uncertainties.

Since the losses of PINNs and APNNs are integrals, we approximate them by the Monte Carlo method by selecting  small number of sub-domains randomly   (batch size) and compute the operator $\left \langle \cdot \right \rangle$ and $\Pi(\cdot)$ with quadrature rule. The optimization problem is solved by the Adam version of the gradient descent method~\cite{kingma2015adam}. {The initial value of parameters set $\theta$ in all numerical experiments are generated by Xavier initialization. All the hyper-parameters are chosen for best performance after trying these experiments.}
Specifically,
in most of our experiments, we fix the $x$ domain as $[0, 1]$ and approximate the problem at  randomly selected points $\{(t_i, x_i, v_i)\}$. It should be pointed out that the integral of $v$ for operator $\mathcal{L} / \left \langle \cdot  \right \rangle$ is computed by quadrature rule (the Gauss-Legendre Integration). In detail, assume $\{{w_i, v_i}'\}_{i = 1}^n$ are the nodes and weights,
which $\{{v_i}'\}_{i = 1}^n$ are the roots of the Legendre polynomials with degree $n$ and $\{{w_i}\}_{i = 1}^n$ are determined for accuracy. Then we can approximate an integral, for example, $\int_{-1}^1 f(v) \diff{v}$ by the summation of linear combination of $f( {v_i}'):\sum_{i = 1}^n w_i f( {v_i}')$.

The empirical risk for PINN is as follows
\begin{equation}
    \begin{aligned}
        \mathcal{R}^{\varepsilon}_{\text{PINN}} = & \frac{1}{N_1} \sum_{i = 1}^{N_1} \left| \varepsilon^2 \partial_t f^{\text{NN}}_{\theta}(t_i, x_i, v_i) + \varepsilon \bm{v} \cdot \nabla_x f^{\text{NN}}_{\theta}(t_i, x_i, v_i) - \mathcal{L} f^{\text{NN}}_{\theta}(t_i, x_i) - \varepsilon^2 Q \right|^2 \\
                                                  & +  \frac{\lambda_1}{N_2}  \sum_{i = 1}^{N_2} |\mathcal{B}f^{\text{NN}}_{\theta}(t_i, x_i, v_i) - F_{\text{B}}(t_i, x_i, v_i)|^2                                                                                                                             \\
                                                  & +  \frac{\lambda_2}{N_3} \sum_{i = 1}^{N_3} |\mathcal{I}f^{\text{NN}}_{\theta}(t_i, x_i, v_i) - f_{0}(t_i, x_i, v_i)|^2,
    \end{aligned}
\end{equation}
where $N_1, N_2, N_3$ are the number of sample points of $\mathcal{T} \times \mathcal{D} \times \Omega, \mathcal{T} \times \partial \mathcal{D} \times \Omega, \mathcal{D} \times \Omega$.

Similarly, the empirical risk for APNN is as follows
\begin{equation}
    \begin{aligned}
        \mathcal{R}^{\varepsilon}_{\text{APNN}} = & \frac{1}{N_1^{(1)}}\sum_{i = 1}^{N_1^{(1)}} | \partial_t \rho^{\text{NN}}_{\theta}(t_i, x_i) + \nabla_x \cdot \left \langle   \bm{v} g^{\text{NN}}_{\theta} \right \rangle (t_i, x_i) - Q |^2    \\
                                                  & +  \frac{1}{N_1^{(2)}} \sum_{i = 1}^{N_1^{(2)}}| \varepsilon^2 \partial_t g^{\text{NN}}_{\theta} (t_i, x_i, v_i) + \varepsilon (I - \Pi)(\bm{v} \cdot \nabla_x g^{\text{NN}}_{\theta})(t_i, x_i) \\
                                                  & \quad  + \bm{v}_i \cdot  \nabla_{\bm{x}} \rho^{\text{NN}}_{\theta}(t_i, x_i) - {\mathcal L} g^{\text{NN}}_{\theta}(t_i, x_i)-(I -\Pi) \varepsilon Q|^2                                           \\
                                                  & +  \frac{\lambda_1}{N_2}  \sum_{i = 1}^{N_2} |\mathcal{B}(\rho^{\text{NN}}_{\theta}(t_i, x_i) + \varepsilon g^{\text{NN}}_{\theta}(t_i, x_i, v_i)) - F_{\text{B}}(t_i, x_i, v_i)|^2              \\
                                                  & +  \frac{\lambda_2}{N_3} \sum_{i = 1}^{N_3} |\mathcal{I}(\rho^{\text{NN}}_{\theta}(t_i, x_i) + \varepsilon g^{\text{NN}}_{\theta}(t_i, x_i, v_i)) - f_{0}(t_i, x_i, v_i)|^2 .
    \end{aligned}
\end{equation}
where $N_1^{(1)}, N_1^{(2)}, N_2, N_3$ are the number of sample points of $\mathcal{T} \times \mathcal{D}, \mathcal{T} \times \mathcal{D} \times \Omega, \mathcal{T} \times \partial \mathcal{D} \times \Omega, \mathcal{D} \times \Omega$.

To investigate the influence of the Monte Carlo method in the integral with respect to $v$, we conduct an extra experiment for the case $\eps = 10^{-8}$ of Problem II. Besides, we show why the conservation is important   for   the same problem.

The reference solutions are obtained by standard finite difference methods. For most of the time we will check the relative $\ell^2$ error of the density $\rho(x)$ between DNN methods and reference solutions, e.g. for $1$d case,
{
        \begin{equation}
            \text{error} := \sqrt{
                \frac{\sum_j |\rho^{\text{NN}}_{\theta, j} - \rho^{\text{ref}}_j|^2}
                {\sum_j |\rho^{\text{ref}}_j|^2}
            }.
        \end{equation}}

\subsection{One-dimensional problems}

We shall consider different problems in slab geometry from the rarefied regimes ($\varepsilon \approx O(1)$) to the diffusive regimes ($\varepsilon \to 0$). Various boundary conditions and initial conditions will be used and both the transient and the steady state solutions will be presented with different $\varepsilon$'s.

\paragraph{Problem I. Smooth initial data with periodic BC}

We start from the rarefied regimes where $\varepsilon=1$ and consider periodic boundary condition with a smooth initial data as follows
\begin{equation}
    f_0(x, v) = \frac{\rho(x)}{\sqrt{2\pi}}e^{-\frac{v^2}{2}},
\end{equation}
where
\begin{equation}
    \rho(x) = 1 + \cos (4 \pi x).
\end{equation}
The source term, scattering and absorbing coefficients are set as
\begin{equation}
    \sigma_S = 1, \quad \sigma_A = 0, \quad Q = 0, \quad \varepsilon = 1.
\end{equation}

{Enforcing exact periodic boundary is applied to improve the numerical performance. The ansatz is based on a Fourier basis and one can construct a transform $T: x \to \{\sin( 2\pi j x) , \cos(2\pi j x)\}_{j = 1}^k$ before the first layer of DNN ~\cite{han2020solving, CSIAM-AM-2-748}. Here, we set $k = 2$.
}

\begin{figure}[ht]
    \centering
    \includegraphics[width=0.45\textwidth]{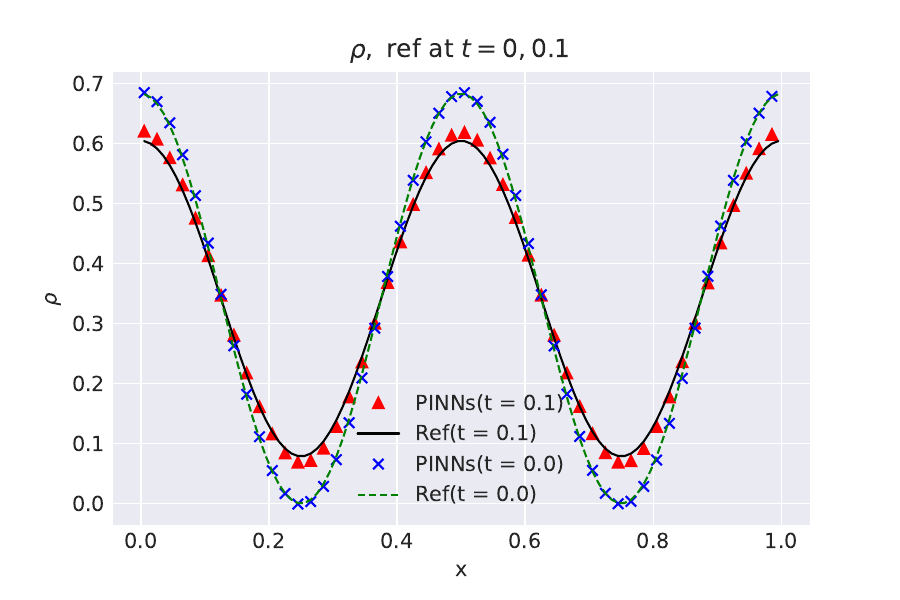}
    \includegraphics[width=0.45\textwidth]{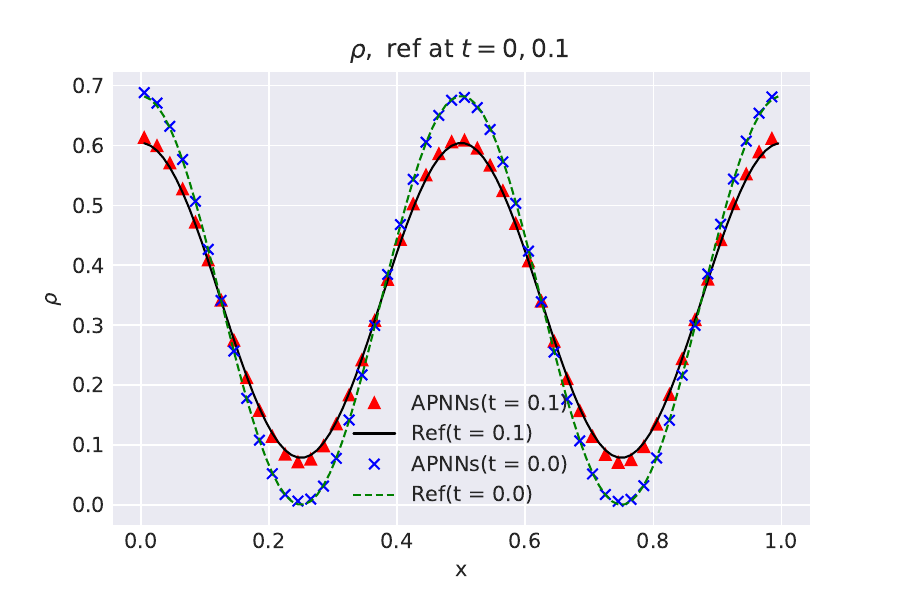}
    \caption{{Problems I. Plot of density $\rho$ at $t = 0.0, 1.0$. Left: PINNs vs. Ref, Right: APNNs vs. Ref. Neural networks are $[5, 128, 128, 128, 128, 1]$ for $\rho$ and $[6, 256, 256, 256, 256, 1]$ for $g, f$. Batch size is $2000$ in domain, $1000$ with penalty $\lambda_2 = 1000$ for initial condition, the number of quadrature points is $30$. Relative $\ell^2$ error of PINNs and APNNs are $2.41 \times 10^{-2}, 1.32  \times 10^{-2}$.}}
    \label{fig:ex1}
\end{figure}

The result is shown in Figure~\ref{fig:ex1} where we can find both PINN and APNN  perform well.

\begin{figure}[htbp]
    \centering
    \subfigure[Density $\rho$ at $t = 1.0$ with $\varepsilon = 1$.]{
        \includegraphics[width=0.4\textwidth]{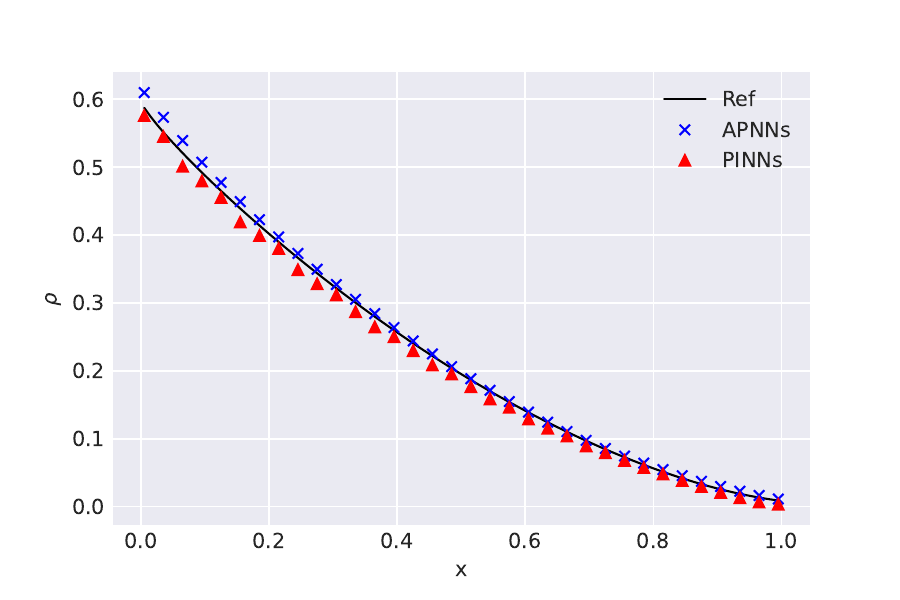}
    }
    \subfigure[Density $\rho$ at $t = 0.5$ with $\varepsilon = 10^{-1}$.]{
        \includegraphics[width=0.4\textwidth]{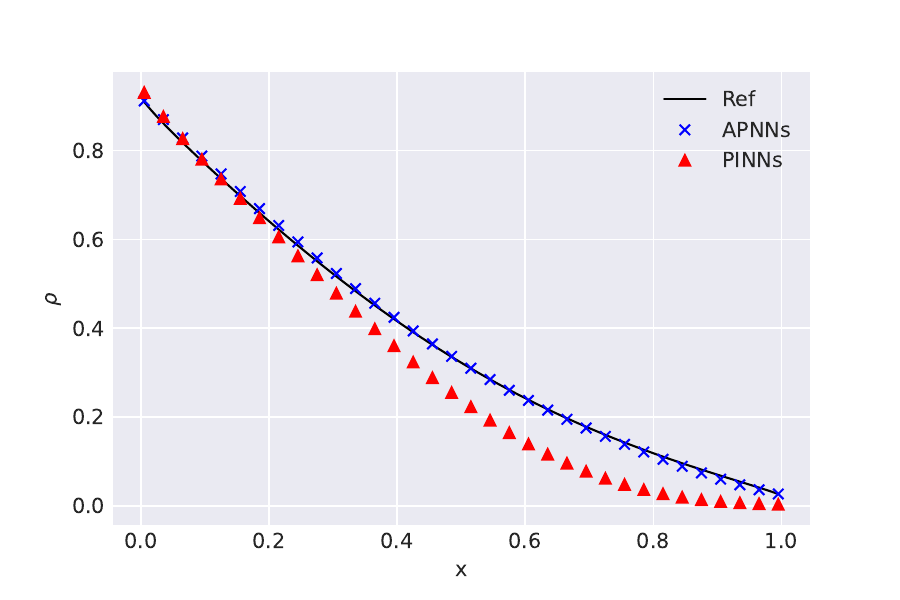}
    }

    \subfigure[Density $\rho$ at $t = 0.1$ with $\varepsilon = 10^{-3}$.]{
        \includegraphics[width=0.4\textwidth]{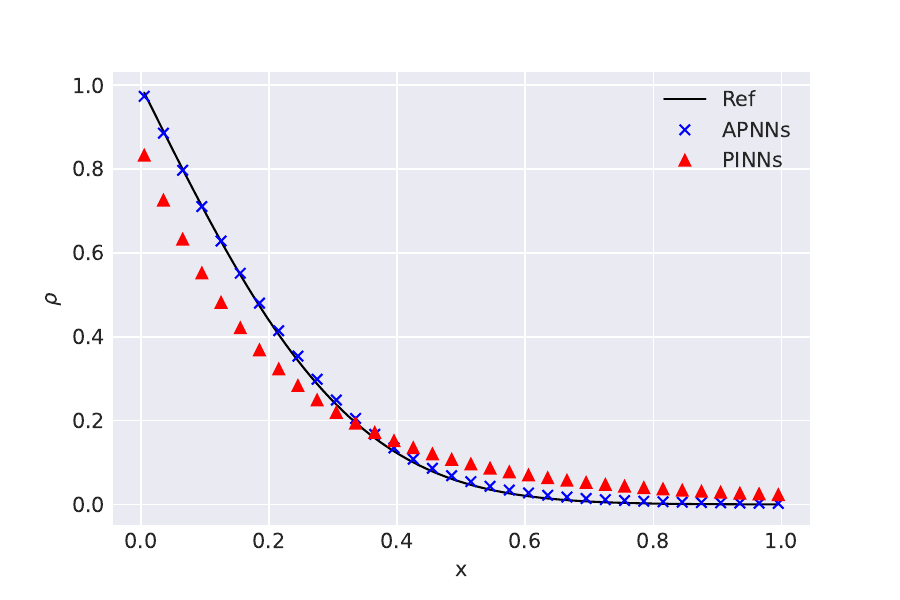}
    }
    \subfigure[Density $\rho$ at $t = 0.05, 0.1$ with $\varepsilon = 10^{-8}$.]{
        \includegraphics[width=0.4\textwidth]{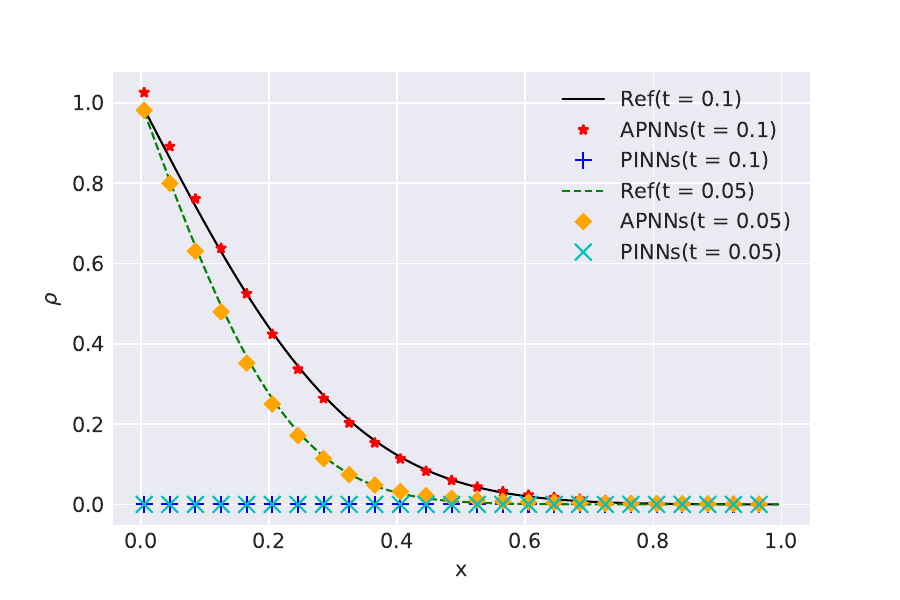}
    }

    \caption{{Problems II. Plot of density $\rho$ for PINNs, APNNs and reference solutions. The number of quadrature points is $30$. Neural networks are ($a$) to ($b$): $[2, 128, 128, 128, 128, 1]$ for $\rho$ and $[3, 256, 256, 256, 256, 1]$ for $g, f$. Batch size is $1000$ for ($a$) to ($d$) in domain, ($a$) and ($b$): $400 \times 2$ with $\lambda_1 = 1$; ($c$) and ($d$): $400 \times 2$ with $\lambda_1 = 10$ for boundary condition, ($a$) and ($b$): $500$ with $\lambda_2 = 1$; ($c$) and ($d$): $1000$ with $\lambda_2 = 10$ for initial condition. Relative $\ell^2$ errors of PINNs and APNNs are ($a$): $4.01 \times 10^{-2}, 1.36  \times 10^{-2}$; ($b$): $1.17 \times 10^{-1}, 3.30  \times 10^{-2}$; ($c$): $2.17 \times 10^{-1}, 1.98  \times 10^{-2}$; ($d$): $9.40 \times 10^{-1}, 2.76  \times 10^{-2}$,
                respectively.}}
    \label{fig:ex2}
\end{figure}

\paragraph{Problem II. In-flow boundary condition}
In the second example we consider the isotropic in-flow boundary conditions:
\begin{equation}
    x \in [0, 1], \quad F_L(v) = 1, \quad F_R(v) = 0.
\end{equation}
{and initial condition $f_0(x, v) = 0$.} The source term, scattering and absorbing coefficients are set as
\begin{equation}
    Q = 0, \quad \sigma_S = 1, \quad \sigma_A = 0, \quad \varepsilon = 1,  10^{-1},  10^{-3},  10^{-8}.
\end{equation}
{The results are shown in Figure \ref{fig:ex2} with exactly $\left \langle g \right \rangle = 0$. Table \ref{tab:ex2} records the relative $\ell^2$ errors of PINNs and APNNs in terms of penalty parameters $(\lambda_1, \; \lambda_2)$ with $\varepsilon = 10^{-3}$. One can find that the approximate accuracy is not good if $\lambda_1, \lambda_2$ are too small or too large.}

{
\begin{table}[tbhp]
    \centering
    \begin{tabular}{ccc}
        \toprule[1pt]
        \noalign{\smallskip}
        \multirow{2}*{$(\lambda_1, \; \lambda_2)$}
         & \multicolumn{2}{c}{Error}                             \\
         & \multicolumn{1}{c}{PINNs} & \multicolumn{1}{c}{APNNs} \\
        \noalign{\smallskip}
        \midrule[1pt]
        \noalign{\smallskip}
        \multirow{1}*{(1, \; 10)}
         & $7.09\times 10^{-1}$      & $5.67\times 10^{-2}$      \\
        \multirow{1}*{(100, \;  10)}
         & $3.57\times 10^{-1}$      & $4.72\times 10^{-2}$      \\
        \multirow{1}*{(10, \; 10)}
         & $2.17\times 10^{-1}$      & $1.98\times 10^{-2}$      \\
        \multirow{1}*{(10, \;  1)}
         & $3.18\times 10^{-1}$      & $1.11\times 10^{-1}$      \\
        \multirow{1}*{(10, \;  100)}
         & $5.52\times 10^{-1}$      & $5.59\times 10^{-2}$      \\
        \noalign{\smallskip}
        \bottomrule[1pt]
    \end{tabular}
    \caption{Problems II. Relative $\ell^2$ errors of PINNs and APNNs in terms of penalty parameters $(\lambda_1, \; \lambda_2)$ with $\varepsilon = 10^{-3}$.}\label{tab:ex2}
\end{table}}

\begin{figure}[htbp]
    \centering
    \subfigure{
        \includegraphics[width=0.4\textwidth]{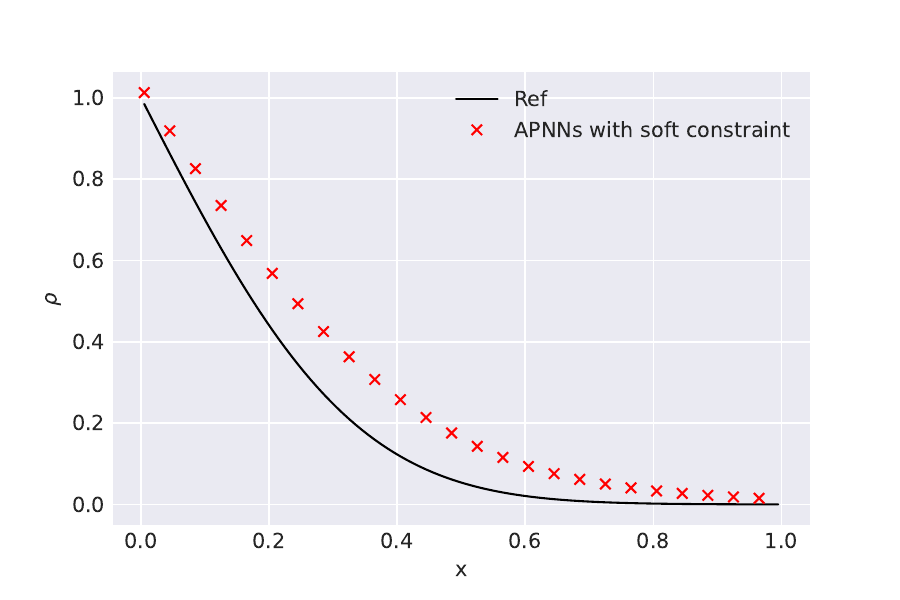}
    }
    \subfigure{
        \includegraphics[width=0.4\textwidth]{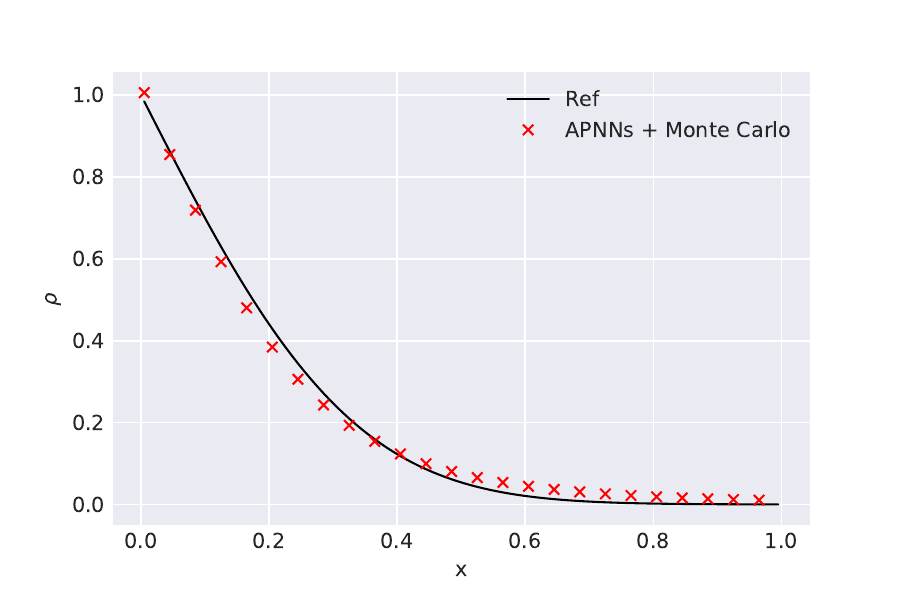}
    }
    \caption{{Problems II. Plot of density $\rho$ at $t = 0.1$. Both neural networks are $[2, 128, 256, 256, 128, 1]$ for $\rho$ and $[3, 128, 256, 512, 256, 128, 1]$ for $g$.
                    {\it Left}  : APNNs with soft constraint $\left \langle g \right \rangle = 0$ vs. Ref.  Batch size is $1000$ in domain, $400 \times 2$ with penalty $\lambda_1 = 10$ for boundary condition, $1000$ with penalty $\lambda_2 = 10$ for initial condition and penalty $\lambda_3 = 1$ for constrain $\left \langle g \right \rangle = 0$. Relative $\ell^2$ error of APNNs with constrain $\left \langle g \right \rangle = 0$ is $2.72 \times 10^{-1}$.
                    {\it Right}:  APNNs with integral by Monte Carlo method vs. Ref. Batch size is $800$ in domain, $500 \times 2$ with penalty $\lambda_1 = 1$ for boundary condition and $500$ with penalty $\lambda_2 = 1$ for initial condition, the number of quadrature points is $100$. The number of sample is 100 for each iteration. The relative $\ell^2$ error of APNNs with constraint $\left \langle g \right \rangle = 0$ is $6.73 \times 10^{-2}$.}}
    \label{fig:ex67}
\end{figure}

{In the left of Figure~\ref{fig:ex67}, this conservation condition is not exactly satisfied by treating it as a soft constraint, i.e., using~\eqref{eq:soft-constraint-g}. Clearly failure to conserve the mass gives poor result. In Figure~\ref{fig:ex2}(d), for small $\varepsilon$ PINN does not give accurate results while APNN gives quite accurate results for all $\varepsilon$ tested. The right of Figure~\ref{fig:ex67} shows the result where the integral with respect to $v$ is computed by the Monte Carlo method.}

% The training process about loss for PINNs and APNNs with case $\varepsilon = 10^{-8}$ are plotted in Figure~\ref{fig:process}, which shows that APNN has a much smaller training error.

% \begin{figure}[ht]
%     \centering
%     \includegraphics[width=0.45\textwidth]{Figure/loss_1e-8_pinns.pdf}
%     \includegraphics[width=0.45\textwidth]{Figure/loss_1e-8_apnns.pdf}
%     \caption{Plot of loss for PINNs and APNNs with case $\varepsilon = 10^{-8}$. Relative $\ell^2$ errors of PINNs and APNNs are $9.40 \times 10^{-1}, 2.76  \times 10^{-3}$ respectively.}
%     \label{fig:process}
% \end{figure}

% \begin{figure}[ht]
%     \centering
%     \includegraphics[width=0.6\textwidth]{Figure/ex_2_1e-8_noexact.pdf}
%     \caption{Problems II. Plot of density $\rho$ at $t = 0.1$: APNNs with soft constraint $\left \langle g \right \rangle = 0$ vs. Ref. $\varepsilon = 10^{-8}$ and neural networks are $[2, 128, 256, 256, 128, 1]$ for $\rho$ and $[3, 128, 256, 512, 256, 128, 1]$ for $g$. Batch size is $1000$ in domain, $400 \times 2$ with penalty $\lambda_1 = 10$ for boundary condition, $1000$ with penalty $\lambda_2 = 10$ for initial condition and penalty $\lambda_3 = 1$ for constrain $\left \langle g \right \rangle = 0$. Relative $\ell^2$ error of APNNs with constrain $\left \langle g \right \rangle = 0$ is $2.72 \times 10^{-1}$.}
%     \label{fig:ex6}
% \end{figure}

% \begin{figure}[ht]
%     \centering
%     \includegraphics[width=0.8\textwidth]{Figure/ex_2_1e-8_mc.pdf}
%     \caption{Problems II. Plot of density $\rho$ at $t = 0.1$: APNNs with integral by Monte Carlo method vs. Ref. $\varepsilon = 10^{-8}$ and neural networks are $[2, 128, 256, 256, 128, 1]$ for $\rho$ and $[3, 128, 256, 512, 256, 128, 1]$ for $g$. Batch size is $800$ in domain, $500 \times 2$ with penalty $\lambda_1 = 1$ for boundary condition and $500$ with penalty $\lambda_2 = 1$ for initial condition, the number of quadrature points is $100$. The number of sample is 100 for each iteration. The relative $\ell^2$ error of APNNs with constraint $\left \langle g \right \rangle = 0$ is $6.73 \times 10^{-2}$.}
%     \label{fig:ex7}
% \end{figure}

\paragraph{Problem III. A variable scattering coefficient}

\begin{figure}[htbp]
    \centering
    \includegraphics[width=0.4\textwidth]{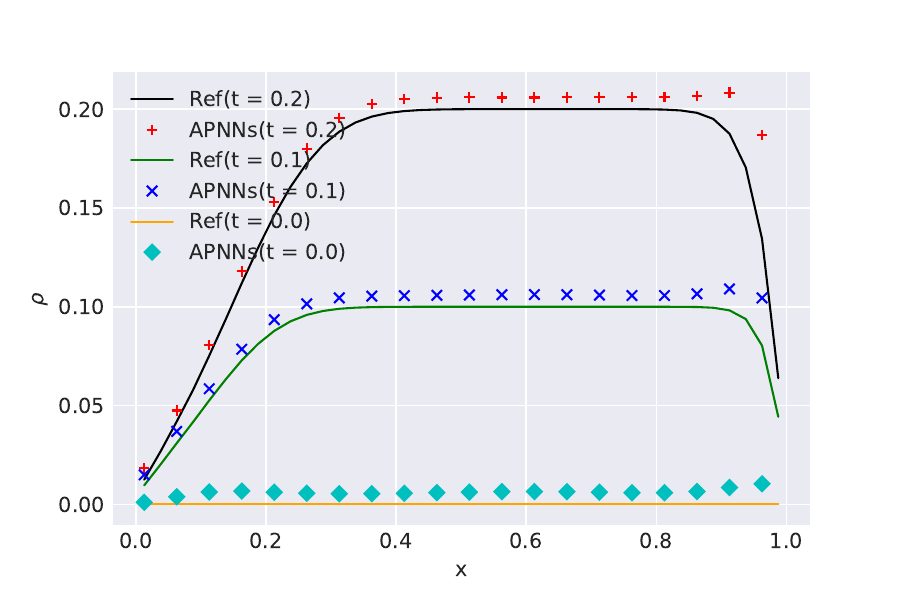}
    \caption{{Problems III. Plot of density $\rho$ with $\varepsilon = 10^{-2}$ at $t = 0.0, 0.1, 0.2$ for APNNs and Reference solutions. Neural networks are $[2, 128, 128, 128, 128, 1]$ for $\rho$ and $[3, 256, 256, 256, 256, 1]$ for $g$. Batch size is $500$ in domain, $200 \times 2$ with penalty $\lambda_1 = 1$ for boundary condition and $200$ with penalty $\lambda_2 = 1$ for initial condition, the number of quadrature points is $30$.}}
    \label{fig:ex8}
\end{figure}

Let
\begin{equation}
    x \in [0, 1], \quad F_L(v) = 1, \quad F_R(v) = 0,
\end{equation}
{and initial condition $f_0(x, v) = 0$.} The source term, scattering and absorbing coefficients are set as
\begin{equation}
    Q = 1, \quad \sigma_S = 1 + (10x)^2, \quad \sigma_A = 0, \quad \varepsilon = 0.01.
\end{equation}
In Figure~\ref{fig:ex8} we report the numerical solution  by APNN at time $t=0.0, 0.1, 0.2$.
In this problem we have a source term and the scattering cross section that depend on $x$, so the scaling term $\sigma_S / \varepsilon$ ranges from $1/\varepsilon \to O(1)$, a problem with mixing scales.{The numerical results how reasonably good performance of APNN.}

\paragraph{Problem IV. A problem with boundary layer}

\begin{figure}[htbp]
    \centering
    \includegraphics[width=0.4\textwidth]{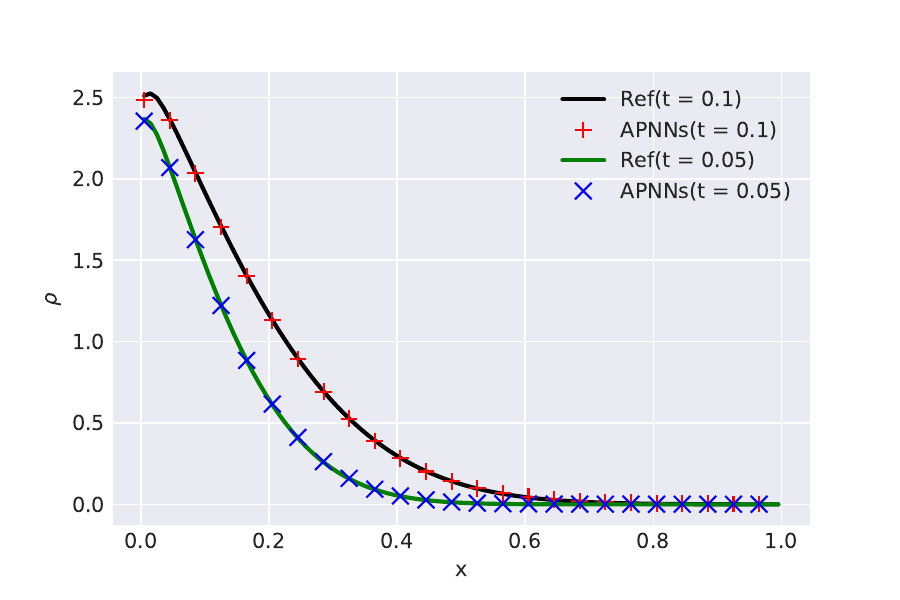}
    \caption{{Problems IV. Plot of density $\rho$ with $\varepsilon = 5 \times 10^{-2}$ at $t = 0.05, 0.1$ for APNNs and Reference solutions.. Neural networks are $[2, 128, 128, 128, 128, 1]$ for $\rho$ and $[3, 256, 256, 256, 256, 1]$ for $g$. Batch size is $1000$ in domain, $400 \times 2$ with penalty $\lambda_1 = 10$ for boundary condition and $1000$ with penalty $\lambda_2 = 10$ for initial condition, the number of quadrature points is $30$. Relative $\ell^2$ error of APNNs is $4.80 \times 10^{-3}$.}}
    \label{fig:ex9}
\end{figure}

Let
\begin{equation}
    x \in [0, 1], \quad F_L(v) = 5 \sin(v), \quad F_R(v) = 0,
\end{equation}
{and initial condition $f_0(x, v) = 0$.} The source term, scattering and absorbing coefficients are set as
\begin{equation}
    Q = 0, \quad \sigma_S = 1, \quad \sigma_A = 0, \quad \varepsilon = 0.05.
\end{equation}
{Here since $F_L$ depends on $v$, there is a boundary layer near $x=0$. Figure \ref{fig:ex9} shows that the boundary layer is well captured by APNN.}

\subsection{Two-dimensional problems}

\paragraph{Problem V. Rarefied regime}

Consider a two dimensional problem with
\begin{equation}
    \Gamma =  [0, 1] \times [0, 1], \quad F_B(\bm{x}, \bm{v}) = 0, \quad \bm{n} \cdot \bm{v} < 0, \quad \bm{x} \in \partial \Gamma,
\end{equation}
{and initial condition $f_0(x, v) = 0$.}
The source term, scattering and absorbing coefficients are set as
\begin{equation}
    Q = 1, \quad \sigma_S = 1, \quad \sigma_A = 0, \quad  \varepsilon = 1.
\end{equation}
Here $\bm{n}$ denotes the exterior unit normal vector on $\partial \Gamma$.

\begin{figure}[ht]
    \centering
    \centering\includegraphics[width=0.45\textwidth]{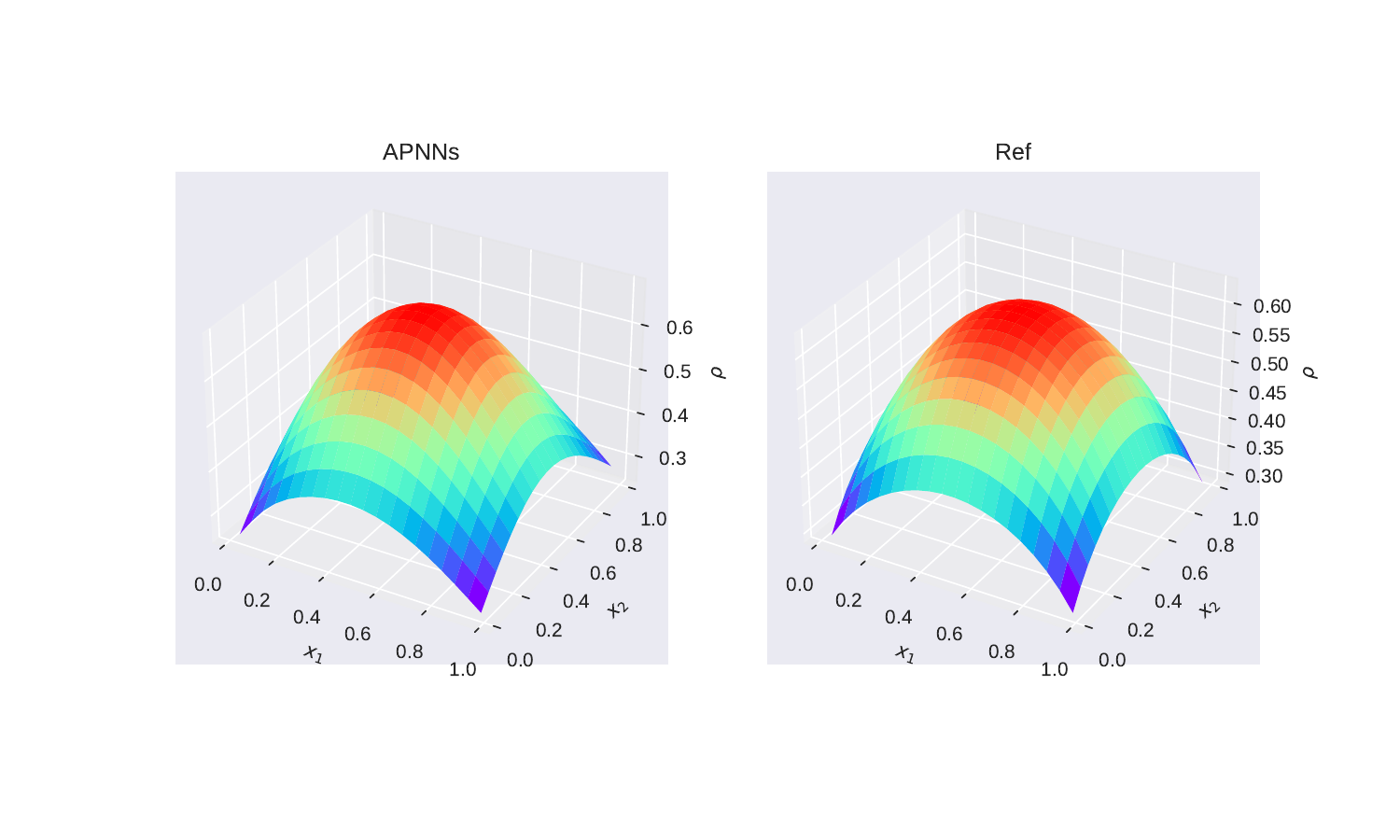}
    \centering\includegraphics[width=0.45\textwidth]{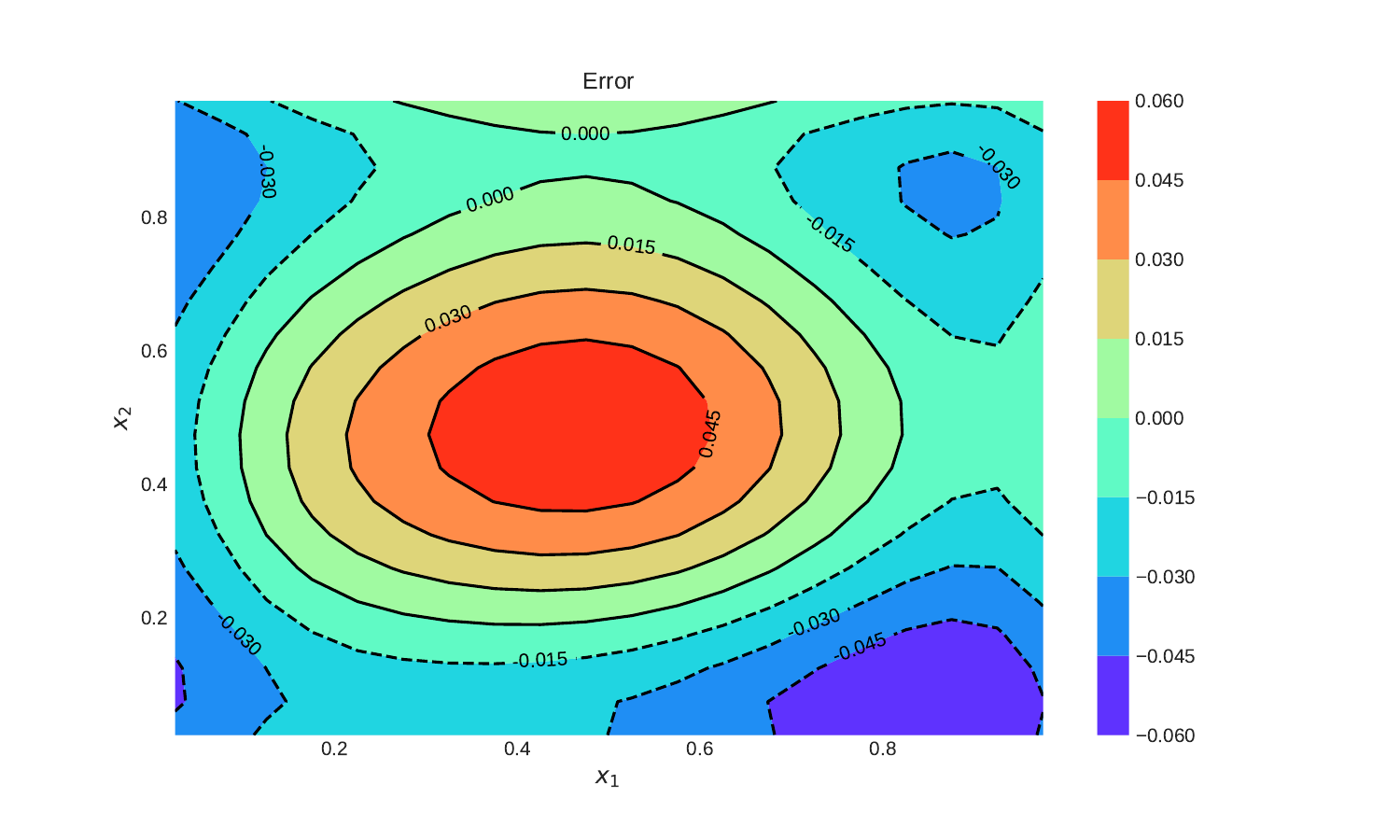}

    \caption{{Problems V. Plot of density $\rho$ with $\varepsilon = 1$ at $t = 1.0$ for APNNs and Reference solution. Neural networks are $[3, 64, 128, 256, 1]$ for $\rho$ and $[5, 64, 128, 256, 512, 1]$ for $g$. Batch size is $800$ in domain, $400 \times 2$ with penalty $\lambda_1 = 1$ for boundary condition and $400$ with penalty $\lambda_2 = 1$ for initial condition, the number of quadrature points is $50$. Relative $\ell^2$ error of APNNs is $5.78 \times 10^{-2}$.}}
    \label{fig:ex10}
\end{figure}

    {Figure~\ref{fig:ex10} shows the density $\rho$ trained by APNNs and reference solution at time $t = 1.0$ and the residual between them.} The results show the good performance of APNNs with absolute error about $0.06$ and relative $\ell^2$ error $5.78 \times 10^{-2}$.

% \begin{figure}[ht]
%     \centering\includegraphics[width=\textwidth]{Figure/ex_5_1e0.pdf}
%     \caption{Problems V. Plot of density $\rho$ at $t = 1.0$: APNNs vs. Ref. $\varepsilon = 1$ and neural networks are $[3, 64, 128, 256, 1]$ for $\rho$ and $[5, 64, 128, 256, 512, 1]$ for $g$. Batch size is $800$ in domain, $400 \times 2$ with penalty $\lambda_1 = 1$ for boundary condition and $400$ with penalty $\lambda_2 = 1$ for initial condition, the number of quadrature points is $50$. Relative $\ell^2$ error of APNNs is $5.78 \times 10^{-2}$.}
%     \label{fig:ex10}
% \end{figure}
% \begin{figure}[ht]
%     \centering\includegraphics[width=\textwidth]{Figure/ex_5_1e0_err.pdf}
%     \caption{Problems VI. Plot of the error of density $\rho$ at $t = 1.0$ between APNNs and Reference.}
%     \label{fig:ex11}
% \end{figure}

\paragraph{Problem VI. A diffusive regime}

\begin{figure}[htbp]
    \centering
    \subfigure[Plot of density for APNNs and Reference solution. {\it Left}  : $Q = 1$; {\it Right}  : $Q = 50$.]{
    \includegraphics[width=0.45\textwidth]{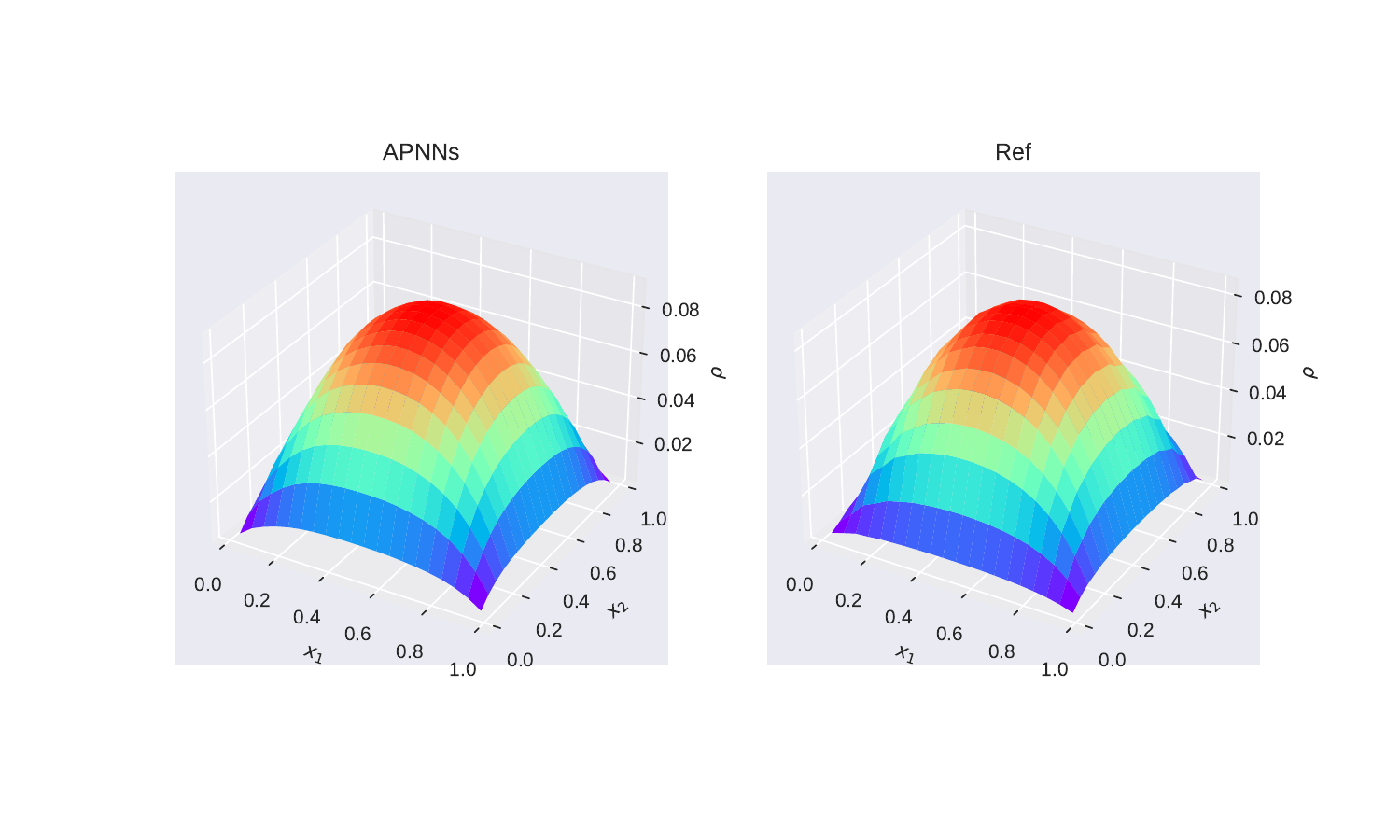}
    \includegraphics[width=0.45\textwidth]{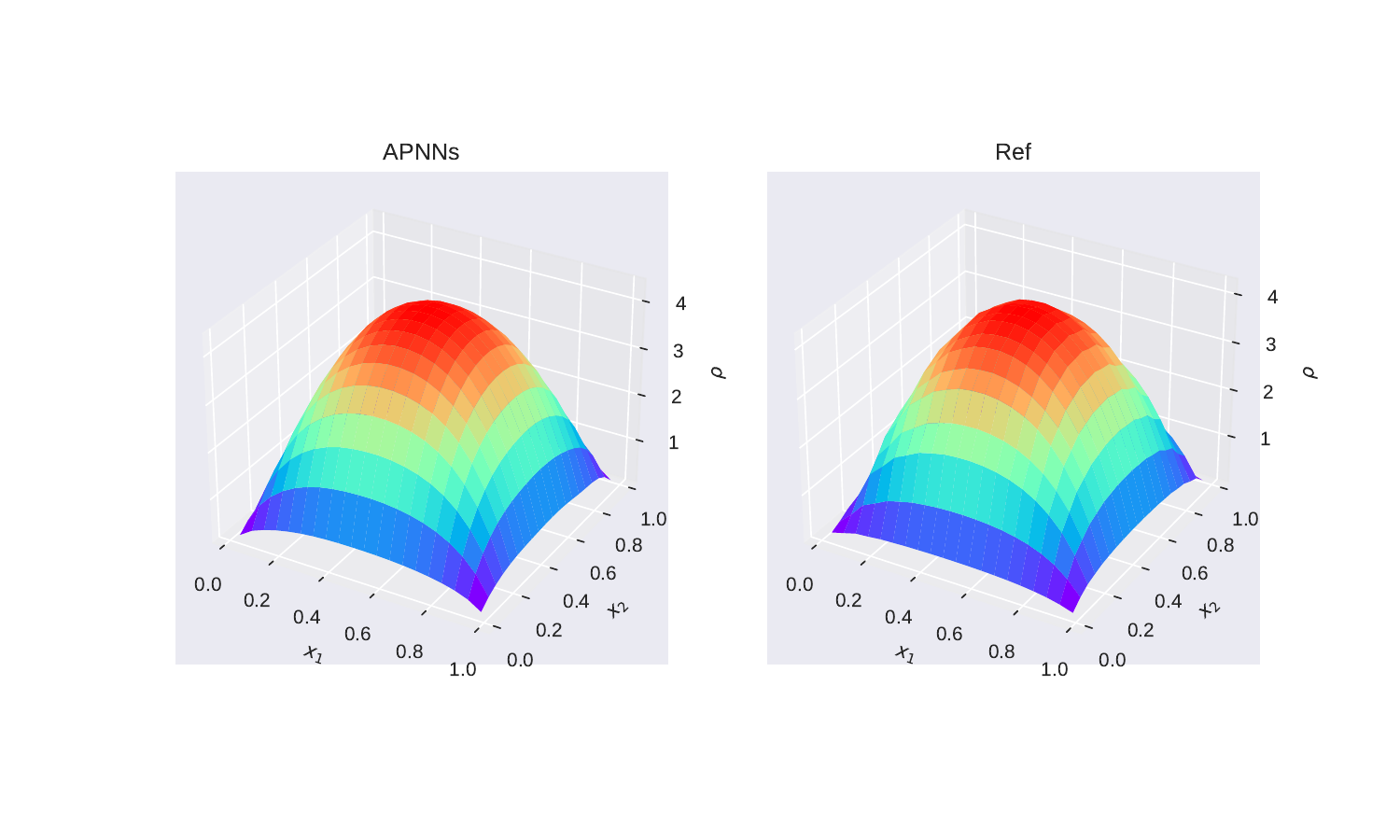}
    }
    \subfigure[Residual. {\it Left}  : $Q = 1$; {\it Right}  : $Q = 50$.]{
    \includegraphics[width=0.45\textwidth]{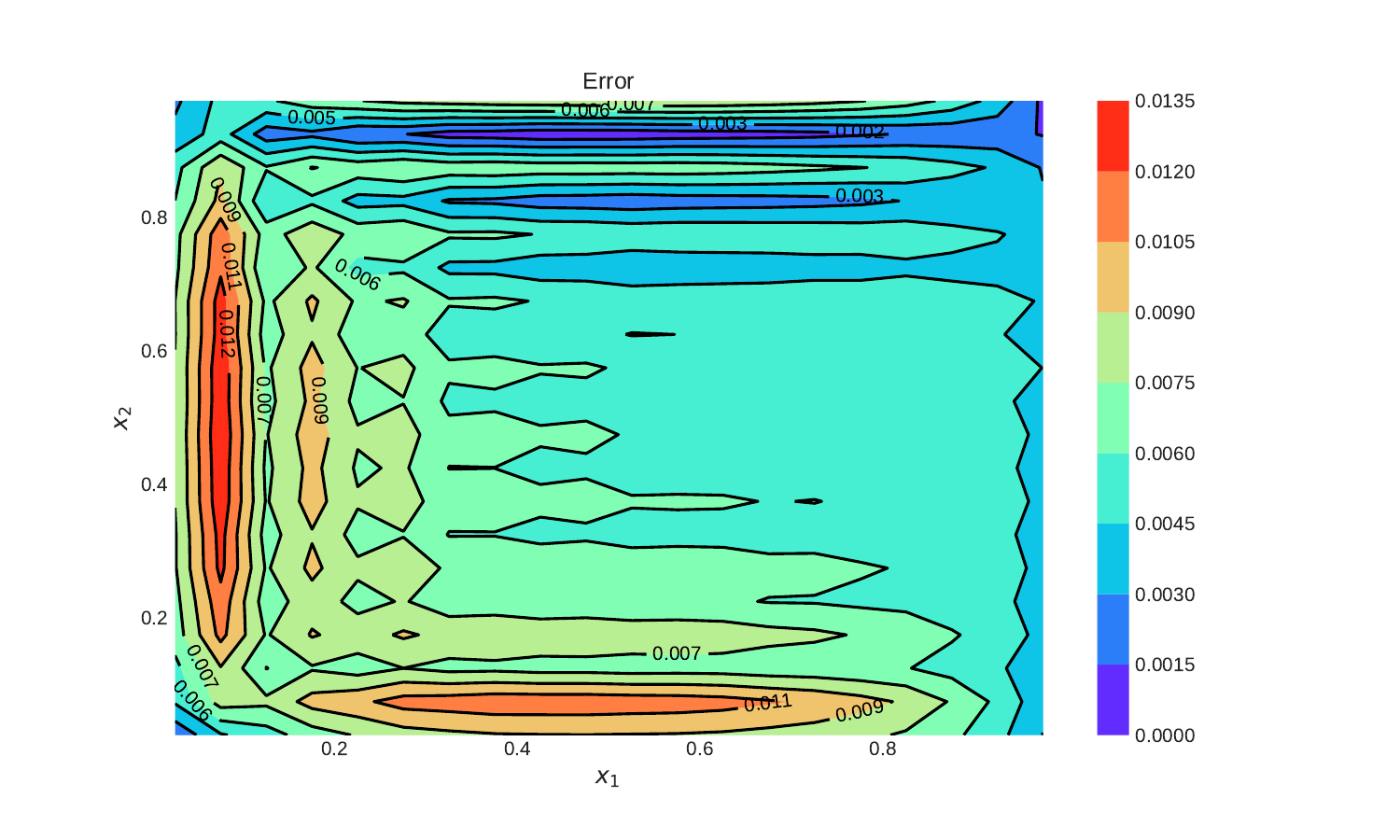}
    \includegraphics[width=0.45\textwidth]{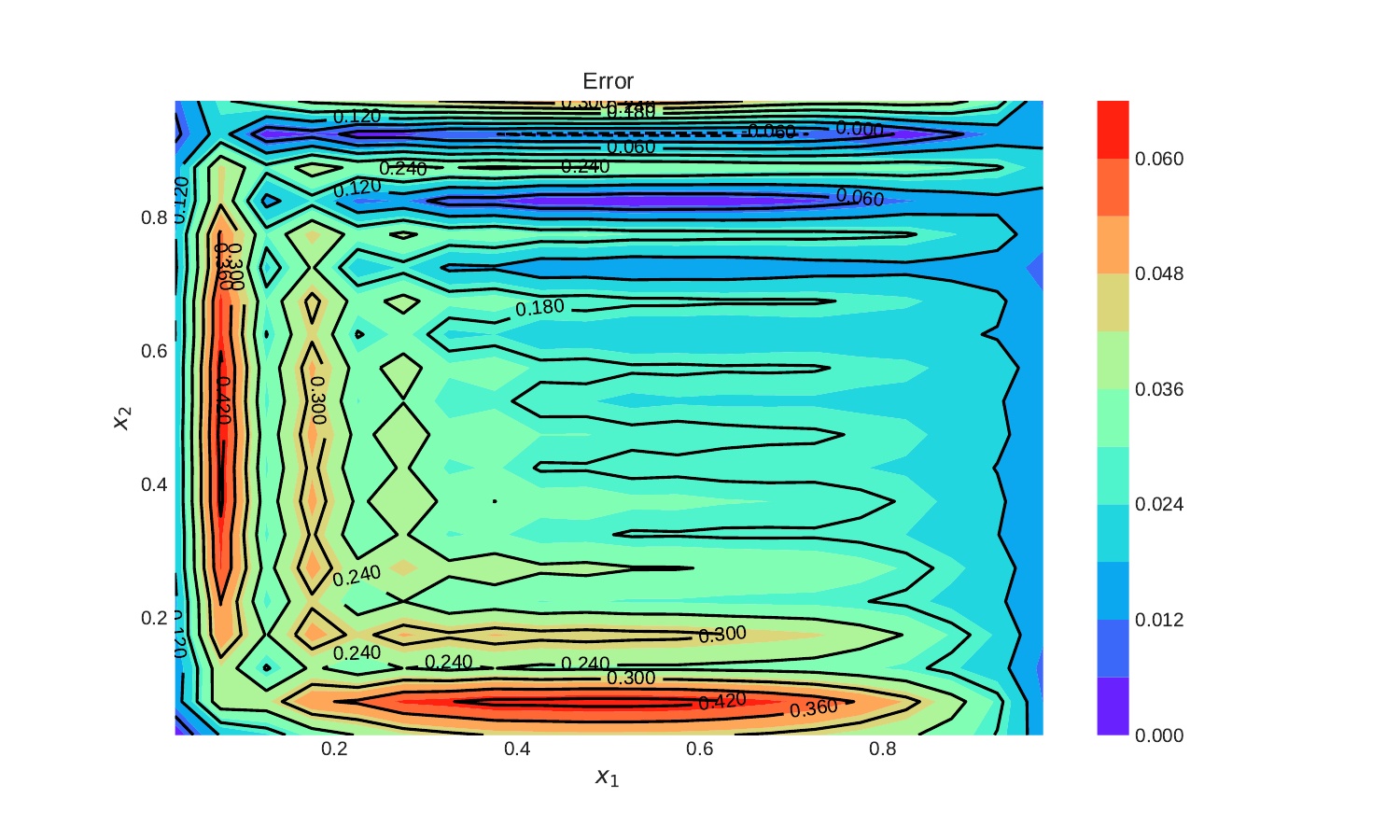}
    }

    \caption{{Problems VI. Plot of density $\rho$ with $\varepsilon = 10^{-8}$ at $t = 0.1$ for APNNs and Reference solution. Neural networks are $[3, 64, 128, 256, 512, 1]$ for $\rho$ and $[5, 64, 128, 256, 512, 1024, 1]$ for $g$. Batch size is $800$ in domain, $400 \times 2$ with penalty $\lambda_1 = 1$ for boundary condition and $400$ with penalty $\lambda_2 = 1$ for initial condition, the number of quadrature points is $50$. Relative $\ell^2$ error of APNNs is $1.29 \times 10^{-2}$ for $Q = 1$ and $4.35 \times 10^{-2}$ for $Q = 50$.}}

    \label{fig:ex12}
\end{figure}

{
    This test is a two dimensional test in the diffusive regime in which most of the setup are the same as previous example except $\varepsilon = 10^{-8}$ and $Q = 1, 50$. Figure~\ref{fig:ex12} show the density $\rho$ trained by APNNs and reference at time $t = 0.1$ and corresponding residual error.}
% \begin{figure}[ht]
%     \centering\includegraphics[width=\textwidth]{Figure/ex_6_1e-8.pdf}
%     \caption{Problems VI. Plot of density $\rho$ at $t = 0.1$: APNNs vs. Ref. $\varepsilon = 10^{-8}$ and neural networks are $[3, 64, 128, 256, 1]$ for $\rho$ and $[5, 64, 128, 256, 512, 1]$ for $g$. Batch size is $800$ in domain, $400 \times 2$ with penalty $\lambda_1 = 1$ for boundary condition and $400$ with penalty $\lambda_2 = 1$ for initial condition, the number of quadrature points is $50$. Relative $\ell^2$ error of APNNs is $1.29 \times 10^{-1}$.}
%     \label{fig:ex12}
% \end{figure}
% \begin{figure}[ht]
%     \centering\includegraphics[width=\textwidth]{Figure/ex_6_1e-8_err.pdf}
%     \caption{Problems VI. Plot of the error of density $\rho$ at $t = 0.1$ between APNNs and Reference.}
%     \label{fig:ex13}
% \end{figure}

\subsection{Uncertainty quantification (UQ) problems}
For the uncertainty quantification problem we consider the linear transport equation with a sine scattering function ($\varepsilon = 1, {10}^{-3}$) or Gaussian scattering function $\sigma_S(\bm z)$ ($\varepsilon = {10}^{-8}$):
\begin{equation}
    \varepsilon \partial_t f + v\partial_x f = \frac{\sigma_S(\bm z)}{\varepsilon}\left ( \frac{1}{2} \int_{-1}^1 f \, dv' - f \right ), \quad x_L < x < x_R, \quad -1 \leq v \leq 1,
\end{equation}
with scattering coefficients
    {
        \begin{equation}
            \sigma_S(\bm z) = 1 + 0.3 \prod_{i=1}^{10} \sin(\pi z_i), \;\bm{z} = (z_1, z_2, \cdots, z_{10}) \sim \mathcal{U}([-1, 1]^{10}),
        \end{equation}}
or
    {
        \begin{equation}
            \sigma_S(\bm z) = 1 + 0.3 \exp{\left(-\frac{|\bm{z}|^2}{2}\right)}, \;\bm{z} = (z_1, z_2, \cdots, z_{20}) \sim \mathcal{U}([-3, 3]^{20}),
        \end{equation}}
    {initial condition $f_0(x, v) = 0$,} and in-flow boundary condition as,
\begin{equation}
    \begin{aligned}
        f(t, x_L = 0, v) & = F_L(v) = 1\quad \text{for} \quad v > 0, \\
        f(t, x_R = 1, v) & = F_R(v) = 0\quad \text{for} \quad v < 0.
    \end{aligned}
\end{equation}
In this problem the 10-dimensional or 20-dimensional vector $\bm{z}$ represents the random input parameters in a typical uncertain problem setup \cite{JXZ}.
Thus, $m_0 = 12, 13$ or $21, 22$ for $\rho$ and $g$ respectively. To compare numerical results, we evaluate $\rho$ at $t = 0.05, 0.1$ by taking expectation on $10^4$ times simulations for $(z_1, \cdots, z_{10})$ or $(z_1, \cdots, z_{20})$. The goal of these examples is
to show the ability of APNNs for high-dimensional problems.

\begin{figure}[htbp]
    \centering\includegraphics[width=0.4\textwidth]{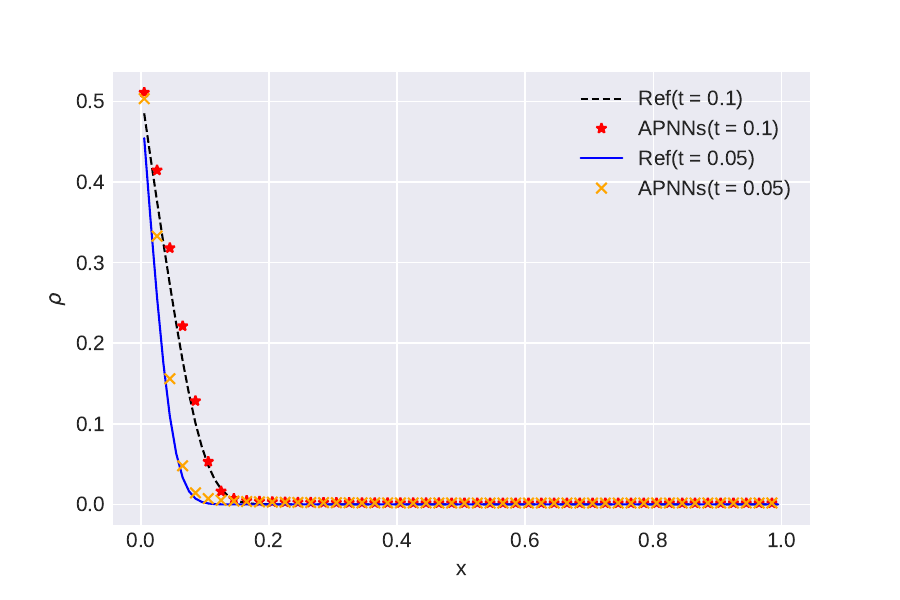}
    \caption{{Problems VII. Plot of density $\rho$ by taking expectation for $\bm{z}$ at $t = 0.05, 0.1$ for APNNs and Reference solution. $\varepsilon = 1, \sigma_S(\bm z) = 1 + 0.3 \prod_{i=1}^{10} \sin(\pi z_i)$ and neural networks are $[12, 128, 256, 256, 512, 1]$ for $\rho$ and $[13, 128, 256, 512, 1024, 1]$ for $g$. Batch size is $1000$ in domain, $500 \times 2$ with penalty $\lambda_1 = 1$ for boundary condition and $500$ with penalty $\lambda_2 = 1$ for initial condition, the number of quadrature points is $30$. Relative $\ell^2$ error of APNNs is $8.78 \times 10^{-2}$.}}
    \label{fig:ex19}
\end{figure}

\paragraph{Problem VII. UQ problem with $\varepsilon=1$}
Let $\varepsilon = 1$ and set scattering coefficients
\begin{equation}
    \sigma_S(\bm z) = 1 + 0.3 \prod_{i=1}^{10} \sin(\pi z_i).
\end{equation}
Figure~\ref{fig:ex19} shows the density $\rho$ trained by APNNs and reference at time $t = 0.05, 0.1$. APNNs gives good approximation results.

\begin{figure}[htbp]
    \centering\includegraphics[width=0.4\textwidth]{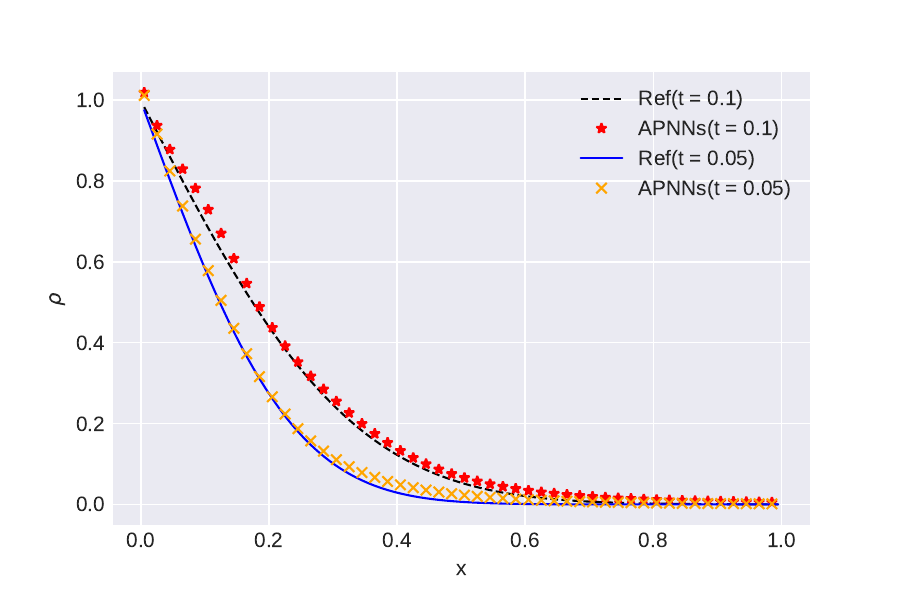}
    \caption{{Problems VIII. Plot of density $\rho$ by taking expectation for $\bm{z}$ at $t = 0.05, 0.1$ for APNNs and Reference solution. $\varepsilon = 10^{-3}, \sigma_S(\bm z) = 1 + 0.3 \prod_{i=1}^{10} \sin(\pi z_i)$ and neural networks are $[12,64, 128, 256, 512, 1]$ for $\rho$ and $[13, 64, 128, 256, 512, 1024, 1]$ for $g$. Batch size is $1000$ in domain, $400 \times 2$ with penalty $\lambda_1 = 10$ for boundary condition and $1000$ with penalty $\lambda_2 = 10$ for initial condition, the number of quadrature points is $30$. Relative $\ell^2$ error of APNNs is $2.75 \times 10^{-2}$. }}
    \label{fig:ex18}
\end{figure}

\paragraph{Problem VIII. UQ problem with $\varepsilon = 10^{-3}$}
Let $\varepsilon = 10^{-3}$ and set scattering coefficients same as Problems VIII.
Figure~\ref{fig:ex18} shows the density $\rho$ trained by APNNs and reference at time $t = 0.05, 0.1$. APNN again gives quite good results.

% \paragraph{Problem X. UQ problem with $\varepsilon = 1$}
% Let $\varepsilon = 1$ and set scattering coefficients same as Problems VIII and IX.
% Figure~\ref{fig:ex19} shows the density $\rho$ train by APNNs and reference at time $t = 0.05, 0.1$. APNN gives good results, too.
% \begin{figure}[ht]
%     \centering\includegraphics[width=\textwidth]{Figure/uq5_1e0.pdf}
%     \caption{Problems X. Plot of density $\rho$ by taking expectation for $\bm{z}$ at $t = 0.05, 0.1$: APNNs vs. Ref. $\varepsilon = 1, \sigma_S(\bm z) = 1 + 0.3 \prod_{i=1}^{10} \sin(\pi z_i)$ and neural networks are $[12, 128, 256, 256, 512, 1]$ for $\rho$ and $[13, 128, 256, 512, 1024, 1]$ for $g$. Batch size is $1000$ in domain, $500 \times 2$ with penalty $\lambda_1 = 1$ for boundary condition and $500$ with penalty $\lambda_2 = 1$ for initial condition, the number of quadrature points is $30$.Relative $\ell^2$ error of APNNs is $8.78 \times 10^{-2}$.}
%     \label{fig:ex19}
% \end{figure}

\paragraph{Problem IX. UQ problem with Gaussian scattering function}
Consider the linear transport equation with $\varepsilon = 10^{-8}$ and scattering coefficients
    {
        \begin{equation}
            \sigma_S(\bm z) = 1 + 0.3 \exp{\left(-\frac{|\bm{z}|^2}{2}\right)}, \;\bm{z} = (z_1, z_2, \cdots, z_{20}) \in [-3, 3]^{20}.
        \end{equation}}

\begin{figure}[htbp]
    \centering\includegraphics[width=0.4\textwidth]{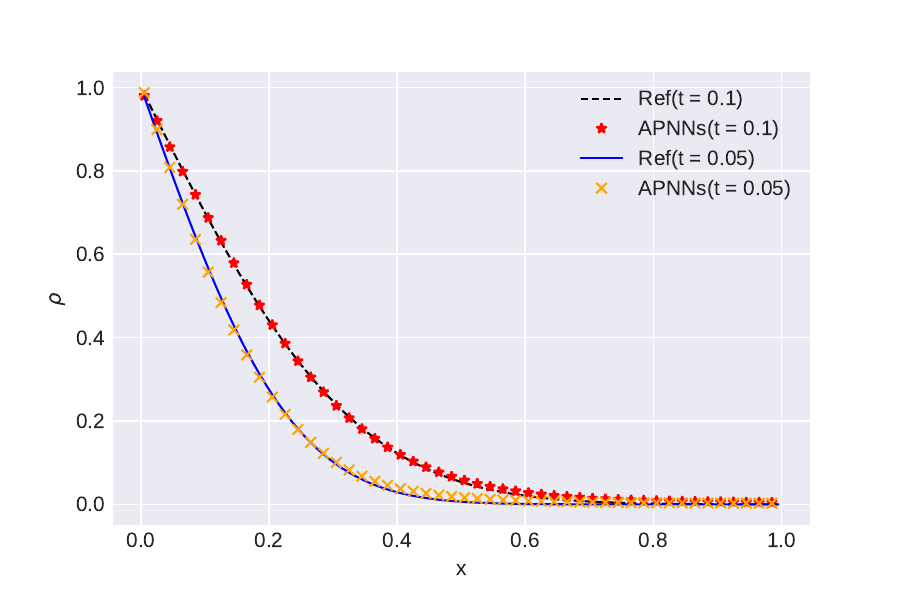}
    \caption{{Problems IX. Plot of density $\rho$ by taking expectation for $\bm{z}$ at $t = 0.05, 0.1$ for APNNs and Reference solution. $\varepsilon = 10^{-8}, \sigma_S(\bm z) = 1 + 0.3 \exp{\left(-{|\bm{z}|^2}/{2}\right)}$ and neural networks are $[22,64, 128, 256, 512, 1]$ for $\rho$ and $[23, 64, 128, 256, 512, 1024, 1]$ for $g$. Batch size is $1000$ in domain, $400 \times 2$ with penalty $\lambda_1 = 10$ for boundary condition and $1000$ with penalty $\lambda_2 = 10$ for initial condition, the number of quadrature points is $30$. Relative $\ell^2$ error of APNNs is $1.51 \times 10^{-2}$.}}
    \label{fig:ex20}
\end{figure}

Figure~\ref{fig:ex20} shows the density $\rho$ trained by APNNs and reference at time $t = 0.05, 0.1$. In this case one can also see good approximation results by APNN.

\section{Conclusion}
In this paper we propose a deep neural network (DNN) for computing the multiscale uncertain linear transport equation with diffusive scaling. Our work follows the framework of Asymptotic-Preserving (AP) schemes for multiscale kinetic equations. We first point out that not all AP schemes will have the desired asymptotic structure
when implementing them in the DNN framework. We design an AP
neural network (APNN) by using the micro-macro decomposition,
together with a mass conservation constraint in the loss function,
that will have the desired AP properties, as will be shown by
various multiscale and high-dimensional uncertain examples.

\section*{Acknowledgement}
This work is partially supported by the National Key R\&D Program of China Project No.~2020YFA0712000 and Shanghai Municipal of Science and Technology Major Project No.~2021SHZDZX0102. Shi Jin is also supported by NSFC grant No.~11871297. Zheng Ma is also supported by NSFC Grant No.~12031013 and partially supported by Institute of Modern Analysis -- A Shanghai Frontier Research Center.

\bibliographystyle{unsrt}
\bibliography{ref.bib}

% \bibliographystyle{plain}
% \bibliography{ref.bib}

\end{document}